\documentclass[aos,nameyear,seceqn,dvips]{arximspdf}

\doi{10.1214/10-AOS860}
\volume{39}
\issue{2}
\pubyear{2011}
\firstpage{887}
\lastpage{930}

\makeatletter
\newtheorem{theorem}{Theorem}
\newtheorem{lemma}{Lemma}
\newtheorem{corollary}{Corollary}
\newtheorem{proposition}{Proposition}
\newproclaim{ree}{Remark}
\newtheorem{cla}{Claim}

\newcommand{\E}{\mathbb{E}}
\newcommand{\N}{\mathbb{N}}
\newcommand{\PP}{\mathbb{P}}
\newcommand{\R}{\mathbb{R}}
\newcommand{\trace}{\operatorname{tr}}
\newcommand{\rank}{\operatorname{rank}}

\newcommand{\AAA}{\mathcal{A}}
\newcommand{\BB}{\mathcal{B}}
\newcommand{\CC}{\mathcal{C}}
\newcommand{\GG}{\mathcal{G}}
\newcommand{\HH}{\mathcal{H}}
\newcommand{\LL}{\mathcal{L}}
\newcommand{\NN}{\mathcal{N}}
\newcommand{\TT}{\mathcal{T}}
\newcommand{\XX}{\mathcal{X}}

\makeatother

\begin{document}
\begin{frontmatter}

\title{Estimation of high-dimensional low-rank matrices\thanksref{T1}}

\thankstext{T1}{Supported in part by the Grant ANR-06-BLAN-0194, ANR ``Parcimonie''
and by the PASCAL-2 Network of Excellence.}
\runtitle{Estimation of high-dimensional low-rank matrices}

\begin{aug}
\author[A]{\fnms{Angelika} \snm{Rohde}\corref{}\ead[label=e1]{angelika.rohde@math.uni-hamburg.de}}
\and
\author[B]{\fnms{Alexandre B.} \snm{Tsybakov}\ead[label=e2]{alexandre.tsybakov@upmc.fr}}

\runauthor{A. Rohde and A. B. Tsybakov}

\affiliation{Universit\"{a}t Hamburg and CREST}
\address[A]{Fachbereich Mathematik\\
Universit\"{a}t Hamburg\\
Bundesstra{\ss}e 55\\
D-20146 Hamburg\\
Germany\\
\printead{e1}} 
\address[B]{Laboratoire de Statistique\\
CREST\\
3, avenue Pierre Larousse\\
F-92240 Malakoff Cedex\\
France\\
\printead{e2}}
\end{aug}

\received{\smonth{12} \syear{2009}}
\revised{\smonth{7} \syear{2010}}

%
\begin{abstract}
Suppose that we observe entries or, more generally, linear
combinations of entries of an unknown $m\times T$-matrix $A$
corrupted by noise. We are particularly interested in the
high-dimensional setting where the number $mT$ of unknown entries
can be much larger
than the sample size $N$. Motivated by several applications, 
we consider estimation of matrix $A$ under the assumption that it has
small rank. This can be viewed as dimension reduction or sparsity
assumption. In order to shrink toward a low-rank representation, we
investigate penalized least squares estimators with a Schatten-$p$
quasi-norm penalty term, $p\leq1$. We study these estimators under two
possible assumptions---a~modified version of the restricted isometry
condition and a uniform bound on the ratio ``empirical norm induced by
the sampling operator/Frobenius norm.'' The main results are stated as
nonasymptotic upper bounds on the prediction risk and on the
Schatten-$q$ risk of the estimators, where $q\in[p,2]$. The rates that
we obtain for the prediction risk are of the form $rm/N$ (for $m=T$),
up to logarithmic factors, where $r$ is the rank of $A$. The particular
examples of multi-task learning and matrix completion are worked out in
detail. The proofs are based on tools from the theory of empirical
processes. As a by-product, we derive bounds for the $k$th entropy
numbers of the quasi-convex Schatten class embeddings
$S_p^M\hookrightarrow S_2^M$, $p<1$, which are of independent interest.
\end{abstract}

%


\end{frontmatter}

\section{Introduction}
\label{sec: introduction}

Consider the observations $(X_i,Y_i)$ satisfying the model
%
\begin{equation}\label{1}
Y_i= \operatorname{tr}(X_i'A^*) +\xi_i,\qquad i=1,\ldots,N,
\end{equation}
where $X_i\in\R^{m\times T}$ are given matrices ($m$ rows, $T$
columns), $A^*\in\R^{m\times T}$ is an unknown matrix, $\xi_i$ are
i.i.d. random errors, $\operatorname{tr}(B)$ denotes the trace of
square matrix $B$ and $X'$ stands for the transposed of $X$. Our aim is
to estimate the matrix $A^*$ and to predict the future $Y$-values based
on the sample $(X_i,Y_i)$, $i=1,\ldots,N$.

We will call model (\ref{1}) the \textit{trace regression model}.
Clearly, for $T=1$ it reduces to the standard regression model. The
``design'' matrices $X_i$ will be called \textit{masks}. This name is
motivated by the fact that we focus on the applications of trace
regression where $X_i$ are very sparse, that is, contain only a small
percentage of nonzero entries. Therefore, multiplication of $A^*$ by
$X_i$ masks most of the entries of $A^*$. The following two examples
are of particular interest.

(i) \textit{Point masks}. For some, typically small, integer $d$
the point masks $X_i$ are defined as elements of the set
\begin{eqnarray*}
&&\XX_d=\Biggl\{\sum_{i=1}^de_{k_i}(m)e_{l_i}'(T)\dvtx1\leq k_i\le
m,1\le l_i\leq T,
\\
&&\hphantom{m,1\le l_i\leq T,tt,}\mbox{with }(k_i,l_i)\not=(k_{i'},
l_{i'}) \mbox{ for }i\not= i'\Biggr\},
\end{eqnarray*}
where $e_k(m)$ are the canonical basis vectors of $\R^m$. In
particular, for $d=1$ the point masks $X_i$ are matrices that have only
one nonzero entry, which equals to~1. The problem of estimation of
$A^*$ in this case becomes the problem of \textit{matrix completion};
the observations $Y_i$ are just some selected entries of $A^*$
corrupted by noise, and the aim is to reconstruct all the entries
of~$A$. The problem of matrix completion dates back at least to Srebro,
Rennie and Jaakkola (\citeyear{SRJbb05a}), Srebro and Shraibman
(\citeyear{SRJgg05b}) and is mainly motivated by applications in
recommendation systems. We will analyze the following two special cases
of matrix completion:
\begin{itemize}
\item[--]\textit{USR (Uniform Sampling at Random) matrix completion.}
The masks $X_i$ are independent, uniformly distributed on
\[
\XX_1=\{e_k(m)e_l'(T)\dvtx 1\leq k\leq m, 1\leq l\leq T\},
\]
and independent from $\xi_1,\dots,\xi_N$.
\item[--]\textit{Collaborative sampling (CS) matrix completion.} The
masks $X_i$ (random or deterministic) belong to $\XX_1$, are all
distinct and independent from $\xi_1,\ldots,\xi_N$.
\end{itemize}
The CS matrix completion model is natural to describe recommendation
systems where every user rates every product only once. The USR matrix
completion can be used for transmission of a large-dimensional matrix
trough a noisy communication channel; only a chosen small number of
entries is transmitted, and nevertheless the original matrix $A^*$ can
be reconstructed by the receiver. An important feature of the
real-world matrix completion problems is that the number of observed
entries is much smaller than the size of the matrix: $N\ll mT$, whereas
$mT$ can be very large. For example, $mT$ is of the order of hundreds
of millions for the Netflix problem.

(ii) \textit{Column or row masks}. If $X_i$ has only a small
number $d$ of nonzero columns or rows, it is called column or row mask,
respectively. We suppose here that $d$ is much smaller than $m$ and
$T$. A remarkable case $d=1$ is covering the problem known in
Statistics and Econometrics as longitudinal (or panel, or
cross-section) data analysis and in Machine Learning as
multi-task learning. In what follows, we will designate this problem as
multi-task learning, to avoid ambiguity. In the simplest version of
multi-task learning, we have $N=nT$ where $T$ is the number of tasks
(for instance, in image detection each task $t$ is associated with a
particular type of visual object, e.g., face, car, chair, etc.), and
$n$ is the number of observations per task. The tasks are characterized
by vectors of parameters $a^*_{t}\in\R^m$, $t=1,\dots,T$, which
constitute the columns of matrix $A^*$:
\[
A^* = (a^*_{1} \cdots a^*_{T}).
\]
The $X_i$ are column masks, each containing only one nonzero column
$\mathbf{x}^{(t,s)}\in\R^m$ (with the convention that $\mathbf{x}^{(t,s)}$ is
the $t$th column):
\[
X_i \in\bigl\{\bigl( 0\cdots 0\underbrace{{\mathbf
x}^{(t,s)}}_{t}0\cdots 0\bigr), t=1,\ldots,T, s=1,\ldots,n \bigr\}.
\]
The column $\mathbf{x}^{(t,s)}$ is interpreted as the vector of predictor
variables corresponding to $s$th observation for the $t$th task. Thus,
for each $i=1,\dots,N$ there exists a pair $(t,s)$ with $t=1,\dots,T,
s=1,\dots,n$, such that
%
\begin{equation}\label{multi}
\trace(X_i'A^*) = (a^*_{t})'\mathbf{x}^{(t,s)}.
\end{equation}
If we denote by $Y^{(t,s)}$ and $\xi^{(t,s)}$ the corresponding values
$Y_i$ and $\xi_i$, then the trace regression model (\ref{1}) can be
written as a collection of $T$ standard regression models:
\[
Y^{(t,s)}=(a^*_{t})'\mathbf{x}^{(t,s)} + \xi^{(t,s)},\qquad t=1,\ldots,T,
s=1,\ldots,n.
\]
This is the usual formulation of the multi-task learning model in the
literature.

For both examples given above, the matrices $X_i$ are sparse in the
sense that they have only a small portion of nonzero entries. On the
other hand, such a sparsity property is not necessarily granted for the
target matrix $A^*$. Nevertheless, we can always characterize $A^*$ by
its rank $r=\rank(A^*)$, and say that a matrix is sparse if it has
small rank; cf. Recht, Fazel and Parrilo (\citeyear{rfp07}). For
example, the problem of estimation of a square matrix
$A^*\in\R^{m\times m}$ is a parametric problem which is formally of
dimension $m^2$ but it has only $(2m-r)r$ free parameters. If $r$ is
small as compared to $m$, then the intrinsic dimension of the problem
is of the order $rm$. In other words, the rank sparsity assumption
$r\ll m$ is a
dimension reduction assumption. 
This assumption will be crucial for the interpretation of our results.
Another sparsity assumption that we will consider is that Schatten-$p$
norm of $A^*$ (see the definition in Section~\ref{sec:prelim} below) is
small for some $0<p\le1$. This is an analog of sparsity expressed in
terms of the $\ell_p$ norm, $0<p\le1$, in vector estimation problems.

Estimation of high-dimensional matrices has been recently studied by
several authors in settings different from the ours [cf., e.g.,
Meinshausen and B\"{u}hlmann (\citeyear{mb06}), Bickel and Levina
(\citeyear{bl08}), Ravikumar et al. (\citeyear{wy85}), Amini and
Wainwright (\citeyear{wainwright09}), Cai, Zhang and Zhou
(\citeyear{caizz10}) and the references cited therein]. Most of
attention was devoted to estimation of a large covariance matrix or its
inverse. In these papers, sparsity is characterized by the number of
nonzero entries of a matrix.

Cand\`{e}s and Recht (\citeyear{cp08}), Cand\`{e}s and Tao
(\citeyear{ct09}), Gross (\citeyear{gff09}), Recht (\citeyear{recht09})
considered the nonnoisy setting ($\xi_i\equiv0$) of the matrix
completion problem under conditions that the singular vectors of $A^*$
are sufficiently spread out on the unit sphere or ``incoherent.'' They
focused on exact recovery of $A^*$. Until now, the sharpest results
are those of Gross (\citeyear{gff09}) and Recht (\citeyear{recht09})
who showed that under ``incoherence condition'' the exact recovery is
possible with high probability if $ N
> Cr(m+T) \log^2 m $ with some constant $C>0$ when we observe $N$
entries of a matrix $A^*\in\R^{m\times T}$ with locations uniformly
sampled at random. Cand\`{e}s and Plan (\citeyear{cp09a}), Keshavan,
Montanari and Oh (\citeyear{kmoff09}) explored the same setting in the
presence of noise, proposed estimators $\hat A$ of $A^*$ and evaluated
their Frobenius norm $\|\hat A- A^*\|_F$. The better bounds are in
Keshavan, Montanari and Oh (\citeyear{kmoff09}) who suggest $\hat A$
such that for $A^*\in\R^{m\times T}$ and $T=\alpha m$ with $\alpha>1$
the squared error $\|\hat A- A^*\|_F^2$ is of the order
$\alpha^{5/2}rm^3(\log N)/N$ with
probability close to 1 when the noise is i.i.d. Gaussian. 

In this paper, we\vspace*{1pt} consider the general noisy setting of the trace
regression problem. We study a class of Schatten-$p$ estimators $\hat
A$, that is, the penalized least squares estimators with a penalty
proportional to Schatten-$p$ norm; cf. (\ref{schest}). The special case
$p=1$ corresponds to the ``matrix Lasso.'' We study the convergence
properties of their prediction error
\[
\hat{d}_{2,N}({\hat A},A^*)^2 =
N^{-1}\sum_{i=1}^N\operatorname{tr}^2\bigl(X_i'({\hat A}-A^*)\bigr)
\]
and of their Schatten-$q$ error. The main contributions of this paper
are the following.
\begin{enumerate}[(a)]
\item[(a)] For all $0<p\le1$, under various assumptions on the masks
$X_i$ (no assumption, USR matrix completion, CS matrix completionmatrix
completionmatrix compl) we obtain different bounds on the prediction
error of Schatten-$p$ estimators involving the Schatten-$p$ norm of
$A^*$.
\item[(b)] For $p$ sufficiently close to 0, under a mild
assumption on $X_i$, we show that Schatten-$p$ estimators achieve the
prediction error rate of convergence $\frac{r\max(m,T)}{N}$, up to a
logarithmic factor. This result is valid for matrices $A^*$ whose
eigenvalues are not exponentially large in $N$. It covers the matrix
completion and high-dimensional multi-task learning problems.
\item[(c)] For all $0<p\le1$, we obtain upper bounds for the
prediction error under the matrix Restricted Isometry (RI) condition on
the masks $X_i$, which is a rather strong condition, and under the
assumption that ${\rm rank}(A^*)\le r$. We also derive the bounds for
the Schatten-$q$ error of ${\hat A}$. The rate in the bounds for the
prediction error is $r\max(m,T)/N$ when the RI condition is satisfied
with scaling factor 1 (i.e., for the case not related to matrix
completion and
high-dimensional multi-task learning). 
%
\item[(d)] We prove the lower bounds showing that the rate
$r\max(m,T)/N$ is minimax optimal for the prediction error and
Schatten-$2$ (i.e., Frobenius) norm estimation error under the RI
condition on the class of matrices $A^*$ of rank smaller than $r$. Our
result is even more general because we prove our lower bound on the
intersection of the Schatten-0 ball with the Schatten-$p$ ball for any
$0<p\le1$, which allows us to show minimax optimality of the upper
bounds of (a) as well. Furthermore, we prove minimax lower bounds for
collaborative sampling and USR matrix completion problems.
\end{enumerate}

The main point of this paper is to show that the suitably tuned
Schatten estimators attain the optimal rate of prediction error up to
logarithmic factors. The striking fact is that we can achieve this not
only under the very restrictive assumption, such as the RI condition,
but also under very mild assumptions on the masks $X_i$.

Finally, it is useful to compare the results for matrix estimation when
the sparsity is expressed by the rank with those for the
high-dimensional vector estimation when the sparsity is expressed by
the number of nonzero components of the vector. For the vector
estimation, we have the linear model
\[
Y_i=X_i'\beta+\xi_i,\qquad i=1,\ldots,N,
\]
where $X_i\in\R^{p}$, $\beta\in\R^p$ and, for example, $\xi_i$ are
i.i.d. $\NN(0,1)$ random variables. Consider the high-dimensional case
$p\gg N$. (This is analogous to the assumption $m^2\gg N$ in the matrix
problem and means that the nominal dimension is much larger than the
sample size.) The sparsity assumption for the vector case has the form
$s\ll N$, where $s$ is the number of nonzero components, or the
\textit{intrinsic dimension} of $\beta$. Let $\hat\beta$ be an
estimator of $\beta$. Then the optimal rate of convergence of the
prediction risk $N^{-1}\sum_{i=1}^N(X_i'(\hat\beta-\beta))^2$ on the
class of vectors $\beta$ with given $s$ is of the order $s/N$, up to
logarithmic factors. This rate is shown to be attained, up to
logarithmic factors, for many estimators, such as the BIC, the Lasso,
the Dantzig selector, Sparse Exponential Weighting, etc.; cf., for
example, Bunea, Tsybakov and Wegkamp (\citeyear{btw07}), Koltchinskii
(\citeyear{kolsf}), Bickel, Ritov and Tsybakov (\citeyear{brt09}),
Dalalyan and Tsybakov (\citeyear{dt098}). Note that this rate is of the
form $\frac{\mathrm{intrinsic\ dimension}}{\mathrm{sample\
size}}=\frac{s}{N}$, up to a logarithmic factor. The general
interpretation is therefore completely analogous to that of the matrix
case: Assume for simplicity that $A^*$ is a~square $m\times m$ matrix
with $\rank(A^*)=r$. As mentioned above, the \textit{intrinsic
dimension} (the number of parameters to be estimated to recover $A^*$)
is then $(2m-r)r$, which is of the order $\sim rm$ if $r\ll m$. An
interesting difference is that the logarithmic risk inflation factor is
inevitable in the vector case [cf. Donoho et~al.
(\citeyear{DonJohHoc92}), Foster and George (\citeyear{fg94})], but not
in the matrix problem, as our results reveal.

This paper is organized as follows. In Section \ref{sec:prelim}, we
introduce notation, some definitions, basic facts about the
Schatten quasi-norms and define the Schatten-$p$ estimators. Section
\ref{sec: Schatten} describes elementary steps in their convergence
analysis and presents two general approaches to upper bounds on the
estimation and prediction error (cf. Theorems~\ref{th:0}
and~\ref{main}) depending on the efficient noise level $\tau$. Our main
results are stated in Sections \ref{sec:bounds_mild_conditions},
\ref{sec:examp} (matrix completion), \ref{sec:multi-task} (multi-task
learning). They are obtained from Theorems~\ref{th:0} and~\ref{main} by
specifying the effective noise level $\tau$ under particular
assumptions on the masks $X_i$. Concentration bounds for certain random
matrices leading to the expressions for the effective noise level are
collected in Section~\ref{sec: stochastic term}. Section \ref{sec:
lower bounds} is devoted to minimax lower bounds. Sections \ref{sec:
prooft1} and \ref{sec:proof_lem} contain the main proofs. Finally, in
Section \ref{sec: embeddings} we establish bounds for the $k$th entropy
numbers of the quasi-convex Schatten class embeddings
$S_p^M\hookrightarrow S_2^M$, $p<1$, which are needed for our proofs
and are of independent interest.

\section{Preliminaries} \label{sec:prelim}

\subsection{Notation, definitions and basic facts}
\label{subsec: notation}

We will write $|\cdot|_2$ for the Euclidean norm in $\R^d$ for any
integer $d$. For any matrix $A\in\R^{m\times T}$, we denote by
$A_{(j,\cdot)}$ for $1\leq j\leq m$ its $j$th row and write
$A_{(\cdot,k)}$ for its $k$th column, $1\leq k\leq T$. We denote by
$\sigma_1(A)\ge\sigma_2(A)\ge\cdots\ge0$ the singular values of $A$.
The (quasi-)norm of some (quasi-) Banach space $\BB$ is canonically
denoted by $\Vert\cdot\Vert_{\BB}$. In particular, for any
matrix $A\in\R^{m\times T}$ and $0<p<\infty$ we consider the Schatten
(quasi-)norms
\[
\Vert A\Vert_{S_p}=\Biggl(\sum_{j=1}^{\min(m,T)}
\sigma_j(A)^{p}\Biggr)^{1/p} \quad\mbox{and}\quad\Vert
A\Vert_{S_{\infty}}=\sigma_1(A).
\]
The Schatten spaces $S_p$ are defined as spaces of all matrices
$A\in\R^{m\times T}$ equipped with quasi-norm $\Vert
A\Vert_{S_p}$. In particular, the Schatten-2 norm coincides with
the Frobenius norm
\[
\Vert A\Vert_{S_2}=\sqrt{\operatorname{tr}(A'A)}
=\biggl(\sum_{i,j} a_{ij}^2\biggr)^{1/2},
\]
where $a_{ij}$ denote the elements of matrix $A\in\R^{m\times T}$.
Recall that for $0<p<1$ the
Schatten spaces $S_p$ are not normed but only quasi-normed, and
$\Vert\cdot\Vert_{S_p}^p$ satisfies the inequality
%
\begin{eqnarray}\label{triang}
\Vert
A+B\Vert_{S_p}^p \leq\Vert A\Vert_{S_p}^p+ \Vert
B\Vert_{S_p}^p
\end{eqnarray}
for any $0<p\le1$ and any two matrices $A,B\in\R^{m\times T}$;
cf. McCarthy (\citeyear{mccarthy67}) and Rotfeld (\citeyear{Rotg69}). 
We will use the following well-known \textit{trace duality} property:
%
\[
|\operatorname{tr}(A'B)| \le\Vert A\Vert_{S_1} \Vert
B\Vert_{S_\infty}\qquad \forall A,B\in\R^{m\times T}.
\]

\subsection{Characteristics of the sampling operator}
\label{subsec:sampl_operator}

Let $\LL\dvtx \R^{m\times T}\rightarrow\R^N$ be the \textit{sampling
operator}, that is, the linear mapping defined by
\[
A\mapsto
(\operatorname{tr}(X_1'A),\ldots,\operatorname{tr}(X_N'A))/\sqrt{N}.
\]
We have
\[
\vert\LL(A)\vert_2^2=
N^{-1}\sum_{i=1}^N\operatorname{tr}^2(X_i'A).
\]
Depending on the context, we also write $\hat{d}_{2,N}(A,B)$ for
$\vert\LL(A-B)\vert_2 $, where $A$ and $B$ are any matrices
in $\R^{m\times T}$. Unless the reverse is explicitly stated, we will
tacitly assume that the matrices $X_i$ are nonrandom.

We will denote by $\phi_{\max}(1)$ the maximal rank-1 restricted
eigenvalue of $\LL$:
%
\begin{equation}\label{eq:phi_max}
\phi_{\max}(1) = \sup_{A\in\R^{m\times T}\dvtx
{\operatorname{rank}}(A)=1}\frac{\vert\LL(A)\vert_2}{\Vert
A\Vert_{S_2}}.
\end{equation}
%
We now introduce two basic assumptions on the sampling operator that
will be used in the sequel. The sampling operator $\LL$ will be called
\textit{uniformly bounded} if there exists a constant $c_0<\infty$ such
that
%
\begin{equation}\label{eq: UB}
\sup_{A\in\R^{m\times
T}\setminus\{0\}}\frac{\vert\LL(A)\vert_2^2}{\Vert
A\Vert_{S_2}^2}\leq c_0\qquad\mbox{uniformly in $m$, $T$ and
$N$.}
\end{equation}
Clearly, if $\LL$ is uniformly bounded, then $\phi_{\max}^2(1) \le
c_0$. Condition (\ref{eq: UB}) is trivially satisfied with $c_0=1$ for
USR matrix completion and with $c_0=1/N$ for CS matrix completion.

The sampling operator $\LL$ is said to satisfy the \textit{Restricted
Isometry condition} RI~($r$,$\nu$) for some integer $1\le r\le
\min(m,T)$ and some $0<\nu<\infty$ if there exists a constant
$\delta_r\in(0,1)$ such that
%
\begin{equation}\label{eq:RIcondition}
(1-\delta_r)\Vert A\Vert_{S_2}\leq\nu\vert
\LL(A)\vert_2\leq(1+\delta_r)\Vert A\Vert_{S_2}
\end{equation}
for all matrices $A\in\R^{m\times T}$ of rank at most $r$.

A difference of this condition from the Restricted Isometry condition,
introduced by Cand\`{e}s and Tao (\citeyear{ct05}) in the vector case
or from its analog for the matrix case suggested by Recht, Fazel and
Parrilo (\citeyear{rfp07}), is that we state it with a \textit{scaling
factor} $\nu$. This factor is introduced to account for the fact that
the masks $X_i$ are typically very sparse, so that they do not induce
isometries with coefficient close to one. Indeed, $\nu$ will be large
in the
examples that we consider below. 

\subsection{Least squares estimators with Schatten penalty} \label{sec:main}

In this\vspace*{1pt} paper, we study the estimators $\hat{A}$ defined as a solution
of the minimization problem
%
\begin{equation}\label{schest}
\min_{A\in\R^{m\times
T}}\Biggl(\frac{1}{N}\sum_{i=1}^N\bigl(Y_i-\operatorname{tr}(X_i'A)\bigr)^2
+ \lambda\Vert A\Vert_{S_p}^p\Biggr)
\end{equation}
with some fixed $0<p\leq1$ and $\lambda>0$. 
The case $p=1$ (matrix Lasso) is of outstanding interest since the
minimization problem is then convex and thus can be efficiently solved
in polynomial time. We call $\hat{A}$ the \textit{Schatten-$p$
estimator}. Such estimators have been recently considered by many
authors motivated by applications to multi-task learning and
recommendation systems. Probably, the first study is due to Srebro,
Rennie and Jaakkola (\citeyear{SRJbb05a}) who dealt with binary
classification and considered the Schatten-1 estimator with the hinge
loss rather than squared loss. Argyriou et al. (\citeyear{argyr07}),
Argyriou, Evgeniou and Pontil (\citeyear{argyr08}), Argyriou, Micchelli
and Pontil (\citeyear{argyr09}), Bach (\citeyear{bach08}), Abernethy et
al. (\citeyear{bachetal09}) discussed connections of (\ref{schest}) to
other related minimization problems, along with characterizations of
the solutions and computational issues, mainly focusing on the convex
case $p=1$. Also for the nonconvex case ($0<p<1$), Argyriou et al.
(\citeyear{argyr07}), Argyriou, Evgeniou and Pontil
(\citeyear{argyr08}) suggested an algorithm of approximate computation
of Schatten-$p$ estimator or its analogs. However, for $0<p<1$ the
methods can find only a local minimum in (\ref{schest}), so that
Schatten estimators with such $p$ remain for the moment mainly of
theoretical value. In particular, analyzing these estimators reveals,
which rates of convergence can, in principle, be attained.

The statistical properties of Schatten estimators are not yet well
understood. To our knowledge, the only previous study is that of Bach
(\citeyear{bach08}) showing that for $p=1$, under some condition on
$X_i'$'s [analogous to strong irrepresentability condition in the vector
case; cf. Meinshausen and B\"{u}hlmann (\citeyear{mb06}), Zhao and Yu
(\citeyear{zy05})], $\operatorname{rank}(A^*)$ is consistently recovered by
$\operatorname{rank}(\hat{A})$ when $m,T$ are fixed and $N\to\infty$. Our
results are of a different kind. They are nonasymptotic and meaningful
in the case $mT\gg N > \max(m,T)$. Furthermore, we do not consider the
recovery of the rank, but rather the estimation and prediction
properties of Schatten-$p$ estimators.

After this paper has been submitted, we became aware of interesting
contemporaneous and independent works by Cand\`{e}s and Plan
(\citeyear{cp09b}), Negahban et al.~(\citeyear{NRWYf09}) and Negahban
and Wainwright (\citeyear{nwf09}). Those papers focus on the bounds for
the Schatten-2 (i.e., Frobenius) norm error of the matrix Lasso
estimator under the matrix RI condition. This is related to the
particular instance of our results in item (c) above with $p=1$ and
$q=2$. Their analysis of this case is complementary to ours in several
aspects. Negahban and Wainwright (\citeyear{nwf09}) derive their bound
under the assumption that $X_i$ are matrices with i.i.d. standard
Gaussian elements and $A^*$ belongs to a Schatten-$p'$ ball with $0\le
p'\le1$, which leads to rates different from ours if $p'\ne0$. An
assumption used in this context in Negahban and Wainwright
(\citeyear{nwf09}) is that $N> mT$ (in our notation), which excludes
the high-dimensional case $mT\gg N$ that we are mainly interested in
Cand\`{e}s and Plan (\citeyear{cp09b}) consider approximately low-rank
matrices, explore the closely related matrix Dantzig selector and
provide lower bounds corresponding to a special case of item (d) above.
The results of these papers do not cover the matrix completion and
multi-task learning problems, which are in the main focus of our study.
We also mention a more recent work by Bunea, She and Wegkamp
(\citeyear{bsw10}) dealing with a special case of our model and
analyzing matrix estimators penalized by the rank.

\section{Two schemes of analyzing Schatten estimators}
\label{sec: Schatten}

In this section, we discuss two schemes of proving upper bounds on the
prediction error of $\hat A$. The first bound involves only the
Schatten-$p$ norm of matrix $A^*$. The second involves only the rank
of $A^*$ but needs the RI condition on the sampling operator. 

We start by sketching elementary steps in the convergence analysis of
Schatten-$p$ estimators. By the definition of $\hat{A}$,
\[
\frac{1}{N}\sum_{i=1}^N\bigl(Y_i-\operatorname{tr}(X_i'\hat{A})\bigr)^2
+ \lambda\Vert\hat{A}\Vert_{S_p}^p \le
\frac{1}{N}\sum_{i=1}^N\bigl(Y_i-\operatorname{tr}(X_i'A^*)\bigr)^2 +
\lambda\Vert A^*\Vert_{S_p}^p.
\]
Recalling that $Y_i=\operatorname{tr}(X_i'A^*)+\xi_i$, we can transform
this by a simple algebra to
%
\begin{equation}\label{eq: BI}
\hat{d}_{2,N}(\hat{A},A^*)^2 \leq
\frac{2}{N}\sum_{i=1}^N\xi_i\operatorname{tr}\bigl((\hat
{A}-A^*)'X_i\bigr)+\lambda(\Vert A^*\Vert_{S_p}^p-
\Vert\hat{A}\Vert_{S_p}^p).
\end{equation}
In the sequel, inequality (\ref{eq: BI}) will be referred to as basic
inequality and the random variable $
N^{-1}\sum_{i=1}^N\xi_i\operatorname{tr}((\hat{A}-A^*)'X_i) $ will be
called the \textit{stochastic term}. The core in the analysis of
Schatten-$p$ estimators consists in proving tight bounds for the
right-hand side of the basic inequality (\ref{eq: BI}). For this
purpose, we first need a control of the stochastic term.
Section~\ref{sec: stochastic term} below demonstrates that such a
control strongly depends on the properties of $\LL$, that is, of the
problem at hand.
In summary, Section~\ref{sec: stochastic term} establishes that, under
suitable conditions, for any $0<p\le1$ the stochastic term can be
bounded for all $\delta>0$ with probability close to 1 as follows:
\begin{eqnarray}\label{stoch}
&&\Biggl\vert
\frac{1}{N}\sum_{i=1}^N\xi_i\operatorname{tr}\bigl(X_i'(\hat
{A}-A^*)\bigr)\Biggr\vert \nonumber
\\[-8pt]\\[-8pt]
&&\qquad\leq
\cases{ \tau\Vert\hat{A}-A^*\Vert_{S_1}, &\quad for $p=1$,\cr
\displaystyle\frac{\delta}{2} \hat{d}_{2,N}(\hat{A},A^*)^2 +\tau\delta^{p-1}\Vert\hat{A}-A^*\Vert_{S_p}^p, &\quad for $0<p<1$,
}\nonumber
\end{eqnarray}
where $0<\tau<\infty$ depends on $m,T$ and $N$. The quantity $\tau$
plays a crucial role in this bound. We will call $\tau$ the
\textit{effective noise level}. Exact expressions for $\tau$ under
various assumptions on the sampling operator $\LL$ and on the noise
$\xi_i$ are derived in Section~\ref{sec: stochastic term}. In Table \ref{t1},
we present the values of $\tau$ for three important examples under the
assumption that $\xi_i$ are i.i.d. Gaussian $\mathcal{N}(0,\sigma^2)$
random variables. In the cases listed in Table \ref{t1}, inequality
(\ref{stoch}) holds with probability $1-\varepsilon$, where
$\varepsilon= (1/C)\exp( -C (m+T))$ (first and third example) and
$\varepsilon=(1/C')(\max(m,T)+1)^{-C'}$ (second example) with constants
$C,C'>0$ independent of $N,m,T$.

\begin{table}
\caption{Effective noise level for uniformly bounded $\LL$, USR
and collaborative sampling matrix completion. Here $M=\max(m,T)$, and
the constants $c>0,c(p)>0$ depend only on $\sigma$}\label{t1}
\begin{tabular*}{\textwidth}{@{\extracolsep{\fill}}lcc@{}}
\hline
\textbf{Assumptions on} $\bolds{X_i}$ & \textbf{Assumptions on} $\bolds{N,m,T,p}$ & \textbf{Value of} $\bolds{\tau}$\\
\hline
Uniformly bounded $\mathcal{L}$ & $0<p\le1$ &$c(p)(M/N)^{1-p/2}$\\
USR matrix completion & $p=1$, $(m+T)mT>N$ & $c\min(M/N,(\log M)/\sqrt{N})$\\
CS matrix completion & $p=1$ & $cM^{1/2}/N$\\
\hline
\end{tabular*}
\end{table}

The following two points will be important to understand the subsequent
results:
\begin{itemize}
\item In this paper, we will always choose the regularization parameter
$\lambda$ in the form $\lambda=4\tau$.
\item With this choice of
$\lambda$, the \textit{effective noise level $\tau$ characterizes the
rate of convergence of the Schatten estimator}. The smaller is $\tau$,
the faster is the rate.
\end{itemize}
In particular, the first line in Table \ref{t1} reveals that when
$M=\max(m,T)<N$ the largest $\tau$ corresponds to $p=1$ and it becomes
smaller when $p$ decreases to 0. This suggests that choosing
Schatten-$p$ estimators with $p<1$ and especially $p$ close to 0 might
be advantageous. Note that the assumption of uniform boundedness of
$\LL$ is very mild. For example, it is trivially satisfied with $c_0=1$
for USR matrix completion and with $c_0=1/N$ for CS matrix completion.
However, in these two cases a specific analysis leads to sharper bounds
on the effective noise level (i.e., on the rate of convergence of the
estimators); cf. the second and third lines of Table~\ref{t1}.

In this section, we provide two bounds on the prediction error of $\hat
A$ with a general effective noise level $\tau$. We then detail them in
Sections \ref{sec:bounds_mild_conditions}--\ref{sec:multi-task} for particular values of $\tau$ depending on the
assumptions on the $X_i$. The first bound involves the Schatten-$p$
norm of matrix $A^*$.

\begin{theorem}
\label{th:0}
Let $A^*\in\R^{m\times T}$, and let $0<p\leq1$. Assume
that (\ref{stoch}) holds with probability at least $1-\varepsilon$ for
some $\varepsilon>0$ and $0<\tau<\infty$. Let $ \hat{A}$ be the
Schatten-$p$ estimator defined as a minimizer of (\ref{schest}) with
$\lambda=4\tau$. Then
%
\begin{equation}\label{7}
\hat{d}_{2,N}(\hat{A},A^*)^2 \leq16\tau\| A^*\|_{S_p}^p
\end{equation}
holds with probability at least $1-\varepsilon$.
\end{theorem}

\begin{pf}
From (\ref{eq: BI}) and (\ref{stoch}) with
$\delta=1/2$ and
$\lambda=4\tau$, we get
\[
\hat{d}_{2,N}(\hat{A},A^*)^2 \leq 8\tau(\|\hat{A}-A^*\|_{S_p}^p
+\|
A^*\|_{S_p}^p- \|
\hat{A}\|_{S_p}^p). 
\]
This and the $p$-norm inequality (\ref{triang}) yield (\ref{7}).
\end{pf}

The bound (\ref{7}) depends on the magnitude of the elements of $A^*$
via $\| A^*\|_{S_p}$. The next theorem shows that under the RI
condition 
this dependence can be avoided, and only the rank of $A^*$ affects the
rate of convergence.

\begin{theorem}
\label{main}
Let $A^*\in\R^{m\times T}$ with $\operatorname{rank}(A^*)\leq
r$, and let $0<p\leq1$. Assume that (\ref{stoch}) holds with
probability at least $1-\varepsilon$ for some $\varepsilon>0$ and
$0<\tau<\infty$. Assume also that the Restricted Isometry condition
RI~$((2+a)r$,$\nu)$ holds with some $0<\nu<\infty$, with a sufficiently
large $a=a(p)$ depending only on $p$ and with $0<\delta_{(2+a)r}\le
\delta_0$ for a sufficiently small $\delta_0=\delta_0(p)$ depending
only on $p$.

Let $ \hat{A}$ be the Schatten-$p$ estimator defined as a minimizer of
(\ref{schest}) with $\lambda=4\tau$.
Then with probability
at least $1-\varepsilon$ we have
%
\begin{eqnarray}
\label{eq:mainn}\hat{d}_{2,N}(\hat{A},A^*)^2
&\leq&
C_1r\tau^{2/(2-p)}\nu^{2p/(2-p)},
\\
\label{eq:mainn1}\Vert\hat{A}-A^*\Vert_{S_q}^q
&\leq&
C_2r\tau^{q/(2-p)}\nu^{2q /(2-p)}\qquad \forall q\in[p,2],
\end{eqnarray}
where $C_1$ and $C_2$ are constants, $C_1$ depends only on $p$ and
$C_2$ depends on $p$ and~$q$.
\end{theorem}

Proof of Theorem \ref{main} is given in Section~\ref{sec: prooft1}. The
values $a=a(p)$ and $\delta_0(p)$ can be deduced from the proof. In
particular, for $p=1$ it is sufficient to take $a=19$.

\begin{ree}\label{re1}
Note that if $\nu=1$ the rates in
(\ref{eq:mainn}) and (\ref{eq:mainn1}) do not depend on $p$ if we
assume in addition the uniform boundedness of $\LL$, which is a very
mild condition. Indeed, taking the value of $\tau$ from the first line
of Table~\ref{t1} we see that $r\tau^{2/(2-p)} \nu^{2p/(2-p)} \sim r M/N$ for
all $0<p\le1$. Thus, under the RI condition, using Schatten-$p$
estimators with $p<1$ does not improve the rate of convergence on the
class of matrices $A^*$ of rank at most $r$.
\end{ree}

\textit{Discussion about the scaling factor $\nu$.}
Remark \ref{re1} deals with the case $\nu=1$, which seems to be not always
appropriate 
for trace regression models. To our knowledge, the only available
examples of matrices $X$ such that the sampling operator $\LL$
satisfies the RI condition with $\nu=1$ are \textit{complete matrices},
that is, matrices with all nonzero entries, which are random and have
specific distributions [typically, i.i.d. Rademacher or Gaussian
entries; cf. Recht, Fazel and Parrilo (\citeyear{rfp07})]. Except for degenerate cases [such as
$N=mT$, the $X_i$ distinct and of the form $\sqrt{N}e_k(m)e_l(T)'$ for
$1\leq k\leq m, 1\leq l\leq T$] the sampling operator $\LL$ defines
typically a restricted isometry with $\nu=1$ only if the matrices $X_i$
contain a considerable number of (uniformly bounded) nonzero entries.

Let us now specify the form of the RI condition in the context of
multi-task learning discussed in the \hyperref[sec: introduction]{Introduction}. Using (\ref{multi})
for a matrix $A=(a_1\cdots a_T)$, we obtain
\begin{eqnarray*}
\vert\LL(A)\vert_2^2
&=&
N^{-1}\sum_{i=1}^N\operatorname{tr}^2(X_i'A)
\\
&=&
N^{-1}\sum_{t=1}^T \sum_{s=1}^n a_{t}'\mathbf{x}^{(t,s)}\bigl(\mathbf{x}^{(t,s)}\bigr)'a_{t}=T^{-1}\sum_{t=1}^T a_{t}' \Psi_t
a_{t},
\end{eqnarray*}
where $\Psi_t= n^{-1}\sum_{s=1}^n \mathbf{x}^{(t,s)}(\mathbf{x}^{(t,s)})'$ is
the Gram matrix of predictors for the $t$th task. These matrices
correspond to $T$ separate regression models. The standard assumption
is that they are normalized so that all the diagonal elements of each
$\Psi_t$ are equal to 1. This suggests that the natural RI scaling
factor $\nu$ for such model is of the order $\nu\sim\sqrt{T}$. For
example, in the simplest case when all the matrices $\Psi_t$ are just
equal to the $m\times m$ identity matrix, we find $ \vert
\LL(A)\vert_2^2 = T^{-1}\sum_{t=1}^T a_{t}' \Psi_t a_{t} =
T^{-1}\|A\|_{S_2}^2. $ Similarly, we get the RI condition with scaling
factor $\nu\sim\sqrt{T}$ when the spectra of all the Gram matrices
$\Psi_t, t=1,\dots,T,$ are included in a fixed interval $[a,b]$ with
$0<a<b<\infty$. However, this excludes the high-dimensional task
regressions, such that the number of parameters $m$ is larger than the
sample size, $m>n$. In conclusion, application of the matrix RI
techniques in multi-task learning is restricted to low-dimensional
regression and the scaling factor is $\nu\sim\sqrt{T}$.

The reason for the failure of the RI approach is that the masks $X_i$ are
sparse. The sparser are $X_i$, the larger is $\nu$. The extreme
situation corresponds to matrix completion problems. Indeed, if $N<mT$,
then there exists a matrix of rank $1$ in the null-space of the
sampling operator $\LL$ and hence the RI condition cannot be satisfied.
For $N\ge mT$ we can have the RI condition with scaling factor $\nu
\sim
\sqrt{mT}$, but $N\ge mT$ means that essentially all the entries are
observed, so that the very problem of completion does not arise.

\section{Upper bounds under mild conditions on the sampling operator}
\label{sec:bounds_mild_conditions}

The above discussion suggests that Theorem \ref{main} and, in general,
the argument based on the restricted isometry or related conditions are
not well adapted for several interesting settings. Motivated by this,
we propose another approach described in the next theorem, which
requires only the comparably mild uniform boundedness condition
(\ref{eq: UB}). For simplicity, we focus on Gaussian errors $\xi_i$.
Set $M=\max(m,T)$.

\begin{theorem}\label{main1}
Let $\xi_1,\dots,\xi_N$ be i.i.d. $\NN(0,\sigma^2)$
random variables. Assume that $M>1$, $N>{\rm{e}}M$ and that the
uniform boundedness condition~(\ref{eq: UB}) is satisfied. Let
$A^*\in\R^{m\times T}$ with $\operatorname{rank}(A^*)\leq r$ and the maximal
singular value $\sigma_1(A^*)\le(N/M)^{C^*}$ for some $0<C^*<\infty$.
Set $p=(\log(N/M))^{-1}$,
$c_{\kappa}=(2\kappa-1)(2\kappa)\kappa^{-1/(2\kappa-1)}$ where
$\kappa=(2-p)/(2-2p)$ and
%
\begin{equation}\label{kappa}
\lambda= 4c_{\kappa}(\vartheta/p)^{1-p/2}\biggl(\frac{M}{N}\biggr)^{1-p/2}
\end{equation}
for some $\vartheta\geq C^2$ and $C$ a universal positive constant
independent of $r$, $M$ and $N$.
Then the Schatten-$p$ estimator $ \hat{A}$ defined as a minimizer of
(\ref{schest}) with $\lambda$ as in (\ref{kappa}) satisfies
%
\begin{equation}\label{eq:main1}
\hat{d}_{2,N}(\hat{A},A^*)^2 \le C_3 \vartheta
\frac{rM}{N}\log\biggl(\frac{N}{M}\biggr)
\end{equation}
with probability
at least $1- C\exp(-\vartheta M/C^2)$ where the positive constant $C_3$
is independent of $r$, $M$ and $N$.
\end{theorem}

\begin{pf}
Inequality (\ref{stoch}) holds
with probability at least $1- C\exp(-\vartheta M/\break C^2)$ by Lemma
\ref{nonconv}. We then use~(\ref{7}) and note that, under our choice of
$p$, $\tau\le c \vartheta M/(Np)$ for some constant $c<\infty$, which
does not depend on $M$ and $N$, and
\[
\|A^*\|_{S_p}^p \le r[\sigma_1(A^*)]^{p} \le r \biggl(\frac{N}{M}\biggr)^{C^*p} =
\exp(C^*)r.
\]
\upqed\end{pf}

Finally, we give the following theorem quantifying the rates of
convergence of the prediction risk in terms of the Schatten norms of
$A^*$. Its proof is straightforward in view of Theorem \ref{th:0} and
Lemmas \ref{lemma:stochterm_matr_compl_1} and \ref{nonconv}.

\begin{theorem}\label{th:00}
Let $\xi_1,\dots,\xi_N$ be i.i.d. $\NN(0,\sigma^2)$ random variables
and $A^*\in\R^{m\times T}$. Then the Schatten-$p$ estimator $ \hat{A}$
has the following properties:
\begin{longlist}[(ii)]
\item[(i)] Let $p=1$, and $\lambda=32\sigma\phi_{\max}(1)\sqrt{(m+T)/N}. $
Then
%
\begin{equation}\label{77}
\hat{d}_{2,N}(\hat{A},A^*)^2 \leq C \sigma\phi_{\max}(1)\| A^*\|_{S_1}
\sqrt{\frac{m+T}{N}}
\end{equation}
with probability at least $1-2\exp\{-(2-\log5)(m+T)\}$ where $C>0$ is
an absolute constant.
\item[(ii)] Let $0<p< 1$ and let the uniform boundedness condition (\ref{eq:
UB}) hold. Set $\lambda$ as in (\ref{kappa}). Then
%
\begin{equation}\label{78}
\hat{d}_{2,N}(\hat{A},A^*)^2 \leq C \| A^*\|_{S_p}^p
\biggl(\frac{M}{N}\biggr)^{1-p/2}
\end{equation}
with probability at least $1- C\exp(-\vartheta M/C^2)$ where
$M=\max(m,T)$ and the constant $C>0$ is independent of $r$, $M$ and
$N$.
\end{longlist}
\end{theorem}

In Theorem \ref{thm: lower bound 1} below we show that these rates are
optimal in a minimax sense on the corresponding Schatten-$p$ balls for
the sampling operators satisfying the RI condition.

\section{Upper bounds for noisy matrix completion}
\label{sec:examp}

As discussed in Section~\ref{sec: Schatten}, for matrix completion
problems the restricted isometry argument as in Theorem~\ref{main} is
not applicable. We will therefore use Theorems \ref{th:0} and
\ref{main1}. First, combining Theorem~\ref{th:0} with
Lemma~\ref{lemma:stochterm_matr_compl} of Section \ref{sec: stochastic
term} we get the following corollary.

\begin{corollary}[(USR matrix completion)]\label{cor:matr_comp1}
\textup{(i)}
Let the i.i.d. zero-mean random variables $\xi_i$ satisfy the Bernstein
condition (\ref{eq: Bernstein}). Assume that\break $mT(m+T)>N$ and consider
the USR matrix completion model. Let $\tau_2$ be given by
(\ref{eq:stochterm0}) with some $D\ge2$. Then the Schatten-$1$
estimator $ \hat{A}$ defined with $\lambda=4\tau_2$ satisfies
%
\begin{equation}\label{eq:matr_comp11}
\hat{d}_{2,N}(\hat{A},A^*)^2 \le16 {\bar C} \|A^*\|_{S_1}
\frac{m+T}{N}
\end{equation}
with probability at least $1-4\exp\{-(2-\log5)(m+T)\}$, where ${\bar
C}=4\sigma\sqrt{10D} + 8HD$.

\textup{(ii)} Let the i.i.d. zero-mean random variables $\xi_i$ satisfy the
light tail condition~(\ref{eq: LT}), and let $\tau_3$ be given by
(\ref{eq:nemirov}) for some $B>0$. Then the Schatten-$1$ estimator $
\hat{A}$ defined with $\lambda=4\tau_3$ satisfies
%
\begin{equation}\label{eq:matr_comp11_3}
\hat{d}_{2,N}(\hat{A},A^*)^2 \le16 \|A^*\|_{S_1}
\sqrt{B}\frac{\sigma\log(\max(m+1,T+1))}{\sqrt{N}}
\end{equation}
with probability at least $1-(1/C){\max(m+1,T+1)}^{-CB}$ for some
constant $C>0$ which does not depend on $m,T$ and $N$.
\end{corollary}

Next, combining Theorem \ref{main1} with Lemma~\ref{nonconv} of Section
\ref{sec: stochastic term} we get the following corollary.

\begin{corollary}[(USR matrix completion, nonconvex penalty)]\label{cor:matr_comp2}
Let $\xi_1,\dots,\break\xi_N$ be i.i.d. $\NN(0,\sigma^2)$
random variables. Assume that $M=\max(m,T)>1$, $N>{\rm e}M$ and
consider the USR matrix completion model. Let $A^*\in\R^{m\times T}$
with $\operatorname{rank}(A^*)\leq r$ and the maximal singular value
$\sigma_1(A^*)\le(N/M)^{C^*}$ for some $0<C^*<\infty$. Set $p=(\log
(N/M))^{-1}$, $c_{\kappa}=(2\kappa-1)(2\kappa)\kappa^{-1/(2\kappa-1)}$
where $\kappa=(2-p)/(2-2p)$ and
\[
\lambda= 4c_{\kappa}(\vartheta/p)^{1-p/2}\biggl(\frac{M}{N}\biggr)^{1-p/2}
\]
for some $\vartheta\geq C^2$ with a universal constant $C>0$,
independent of $r$, $M$ and~$N$.
Then the Schatten-$p$ estimator $ \hat{A}$ defined as a minimizer of
(\ref{schest}) satisfies
%
\begin{equation}\label{eq:matr_comp}
\hat{d}_{2,N}(\hat{A},A^*)^2 \le C_3 \vartheta
\frac{rM}{N}\log\biggl(\frac{N}{M}\biggr)
\end{equation}
with probability
at least $1- C\exp(-\vartheta M/C^2)$, where the positive constant
$C_3$ is also independent of $r$, $M$ and $N$.
\end{corollary}

Note that the bounds of Corollaries \ref{cor:matr_comp1}(i) and
\ref{cor:matr_comp2} achieve the rate $r\max(m,\break T)/N$, up to logarithmic
factors under different conditions on the maximal singular value of
$A^*$. If $\max(m,T) < N<mT$ then the condition in
Corollary~\ref{cor:matr_comp2} does not imply more than a polynomial in
$\max(m,T)$ growth on $\sigma_1(A^*)$, which is a mild assumption. On
the other hand, (\ref{eq:matr_comp11}) requires uniform boundedness of
$\sigma_1(A^*)$ by some constant to achieve the same rate. However, the
estimators of Corollary \ref{cor:matr_comp2} correspond to nonconvex
penalty and are computationally hard.

We now turn to the collaborative sampling matrix completion. The next
corollary follows from combination of Theorem \ref{th:0} with
Lemmas~\ref{lemma:stochterm_matr_compl} and
\ref{lemma:stochterm_matr_compl_2} of Section~\ref{sec: stochastic
term}.

\begin{corollary}[(Collaborative sampling)]\label{cor:matr_comp3}
Consider the problem of~\mbox{matrix} completion with collaborative sampling.
\begin{longlist}[(ii)]
\item[(i)] Let $\xi_1,\ldots,\xi_N$ be i.i.d. $\NN(0,\sigma^2)$ random variables.
Let $\tau_4$ be given by (\ref{eq:stochterm010}) with some $D\ge2$.
Then the Schatten-$1$ estimator $ \hat{A}$ defined with
$\lambda=4\tau_4$ satisfies
%
\begin{equation}\label{eq:matr_comp112}
\hat{d}_{2,N}(\hat{A},A^*)^2 \le16 {\bar C} \|A^*\|_{S_1}
\frac{\sqrt{m+T}}{N}
\end{equation}
with probability at least
$1-2\exp\{-(D-\log5)(m+T)\}$, where ${\bar C}=8\sigma\sqrt{D}$.
\item[(ii)] Let $\xi_1,\ldots,\xi_N$ be i.i.d. zero-mean random variables
satisfying the Bernstein condition (\ref{eq: Bernstein}). Let $\tau_5$
be given by (\ref{eq:stochterm011}) with some $D\ge2$. Then the
Schatten-$1$ estimator $ \hat{A}$ defined with $\lambda=4\tau_5$
satisfies
%
\begin{equation}\label{eq:matr_comp113}
\hat{d}_{2,N}(\hat{A},A^*)^2 \le64 \|A^*\|_{S_1}
\frac{\sigma\sqrt{2D(m+T)} + 2HD(m+T)}{N}
\end{equation}
with probability at least
$1-2\exp\{-(D-\log5)(m+T)\}$.
\end{longlist}
\end{corollary}

\begin{ree}\label{re2}
Using the inequality $\Vert A\Vert
_{S_1}\leq\sqrt{r}\Vert A\Vert_{S_2}$ for matrices $A$ of
rank at most
$r$, we find that that the bound (\ref{eq:matr_comp112}) is minimax
optimal on the class of matrices
\[
\{A\in\R^{m\times T}\dvtx \operatorname{rank}(A)\leq r,\Vert
A\Vert_{S_2}^2 \leq C\sigma^2 r\max(m,T)\}
\]
for some constant $C>0$, if the masks $X_1,\dots,X_N$ fulfill the
dispersion condition of Theorem \ref{thm: lower bound 2} below. It is
further interesting to note that the construction in the proof of the
lower bound in Theorem \ref{thm: lower bound 2} fails if the
restriction is $\Vert A\Vert_{S_2}^2\le\delta^2$ where
$\delta^2$ of smaller order than $r\max(m,T)$.
\end{ree}

\section{Upper bounds for multi-task learning}
\label{sec:multi-task}

For multi-task learning, we can apply both Theorems~\ref{main} and
\ref{main1}. Theorem~\ref{main} imposes a strong assumption on
the masks $X_i$, namely the RI condition. Nevertheless, the advantage
is that Theorem \ref{main} covers the computationally easy case $p=1$.

\begin{corollary}[(Multi-task learning; RI condition)]\label{cor:ml2}
Let $\xi_1,\dots,\xi_N$ be i.i.d. $\NN(0,\sigma^2)$ random variables.
Consider the multi-task learning problem with $\operatorname{rank}(A^*)\leq
r$. Assume that the spectra of the Gram matrices $\Psi_t$ are uniformly
in $t$ bounded from above by a constant $c_1<\infty$. Assume also that
the Restricted Isometry condition RI~($21r$, $\nu$) holds with some
$0<\nu<\infty$ and with $0<\delta_{21r}\le\delta_0$ for a sufficiently
small $\delta_0$. Set
\[
\lambda= 32\sigma\sqrt{\frac{c_1(m+T)}{nT^2}}.
\]
Let $ \hat{A}$ be the Schatten-$1$ estimator with this parameter
$\lambda$. Then with probability at least $1-2\exp\{-(2-\log5)(m+T)\}$
we have
\begin{eqnarray*}
\hat{d}_{2,N}(\hat{A},A^*)^2&\leq&{\bar C}_1 c_1 \sigma^2 r \nu^{2}
\biggl(\frac{m+T}{nT^2}\biggr),\\
\Vert\hat{A}-A^*\Vert_{S_q}^q&\leq&{\bar C}_2 c_1^{q/2}
\sigma^{q} r \nu^{2q} \biggl(\frac{m+T}{nT^2}\biggr)^{q/2}\qquad \forall q\in[1,2],
\end{eqnarray*}
where ${\bar C}_1$ is an absolute constant and ${\bar C}_2$ depends
only on $q$.
\end{corollary}

The proof of Corollary \ref{cor:ml2} is straightforward in view of Theorem
\ref{main}, Lemma~\ref{lemma:stochterm_matr_compl_1} and the fact
that, under the premise of Corollary \ref{cor:ml2}, we have $
\vert\LL(A)\vert_2^2 = T^{-1}\sum_{t=1}^T a_{t}' \Psi_t
a_{t} \leq(c_1/T)\Vert A\Vert_{S_2}^2 $ for all matrices
$A\in\R^{m\times T}$, so that the sampling operator is uniformly
bounded [(\ref{eq: UB}) holds with $c_0=c_1/T$], and thus
$\phi_{\max}(1)\le\sqrt{c_0}\le\sqrt{c_1/T}$.

Taking in the bounds of Corollary \ref{cor:ml2} the natural scaling
factor $\nu\sim\sqrt{T}$, we obtain the following inequalities:
%
\begin{eqnarray}
\label{eq:bochk3}\hat{d}_{2,N}(\hat{A},A^*)^2&\leq&{\tilde C}_1
\frac{r(m+T)}{nT},\\
\label{eq:bochk4}\frac1{T}\Vert\hat{A}-A^*\Vert_{S_2}^2&\leq&{\tilde C}_2
\frac{r(m+T)}{nT},
\end{eqnarray}
where the constants ${\tilde C}_1$ and ${\tilde C}_2$ do not depend
on $m,T$ and $n$. 

A remarkable fact is that the rates in Corollary \ref{cor:ml2} are free
of logarithmic inflation factor. This is one of the differences between
the matrix estimation problems and vector estimation ones, where the
logarithmic risk inflation is inevitable, as first noticed by Donoho
et~al. (\citeyear{DonJohHoc92}), Foster and George (\citeyear{fg94}). For more details about optimal
rates of sparse estimation in the vector case, see Rigollet and
Tsybakov (\citeyear{rigtsy10}).

Since the Group Lasso is a special case of the nuclear norm penalized
minimization on block-diagonal matrices [cf., e.g., Bach (\citeyear{bach08})]
Corollary~\ref{cor:ml2} and the bounds (\ref{eq:bochk3}),
(\ref{eq:bochk4}) imply the corresponding bounds for the Group Lasso
under the low-rank assumption. To note the difference with from the
previous results for the Group Lasso, we consider, for example, those
obtained in multi-task setting by Lounici et al. (\citeyear{lptv09,lptv10}). The main
difference is that the sparsity index $s$ appearing in Lounici et
al. (\citeyear{lptv09,lptv10}) is now replaced by $r$. In Lounici et al. (\citeyear{lptv09,lptv10}),
the columns $a^*_t$ of $A^*$ are supposed to be sparse, with the sets
of nonzero elements of cardinality not more than $s$, whereas here the
sparsity is characterized by the rank $r$ of $A^*$.

Finally, we give the following result based on application of
Theorem~\ref{main1}.

\begin{corollary}[(Multi-task learning; uniformly bounded $\LL$)]\label{cor:ml3}
Let $\xi_1,\dots,\break\xi_N$ be i.i.d. $\NN(0,\sigma^2)$ random variables,
and assume that $n>{\rm e}$. Consider the multi-task learning problem
with $A^*\in\R^{m\times T}$, $\operatorname{rank}(A^*)\leq r$, such that the
maximal singular value $\sigma_1(A^*)\le n^{C^*}$ for some
$0<C^*<\infty$. Assume that the spectra of the Gram matrices $\Psi_t$
are uniformly in $t$ bounded from above by $c_0T$ where $c_0<\infty$ is
a constant. Set $p=(\log n)^{-1}$,
$c_{\kappa}=(2\kappa-1)(2\kappa)\kappa^{-1/(2\kappa-1)}$ where
$\kappa=(2-p)/(2-2p)$ and
\[
\lambda= 4c_{\kappa}(\vartheta/p)^{1-p/2}\biggl(\frac{1}{n}\biggr)^{1-p/2}
\]
for some $\vartheta\geq C^2$ and a universal constant $C>0$,
independent of $r$, $m$ and~$n$. Then the Schatten-$p$ estimator $
\hat{A}$ with this parameter $\lambda$ satisfies
%
\begin{equation}\label{eq:ml3}
\hat{d}_{2,N}(\hat{A},A^*)^2 \le C_3 \vartheta\frac{rM}{nT}\log n
\end{equation}
with probability
at least $1- C\exp(-\vartheta M/C^2)$ where $M=\max(m,T)$, and the
positive constant $C_3$ is independent of $r$, $m$ and $n$.
\end{corollary}

Corollary \ref{cor:ml3} follows from Theorem \ref{main1}. Indeed, it
suffices to remark that, under the premises of Corollary \ref{cor:ml3},
we have $ \vert\LL(A)\vert_2^2 = T^{-1}\sum_{t=1}^T a_{t}'
\Psi_t a_{t} \leq c_0\Vert A\Vert_{S_2}^2 $ for all matrices
$A\in\R^{m\times T}$, so that the sampling operator is uniformly
bounded; cf.~(\ref{eq: UB}).

For $m=T$, we can write (\ref{eq:ml3}) in the form
%
\begin{equation}\label{eq:ml4}
\hat{d}_{2,N}(\hat{A},A^*)^2 \le C_3' \frac{rm}{nT}\log n.
\end{equation}
Clearly, this bound achieves the optimal rate ``intrinsic
dimension/sample size'' \mbox{$\sim rm/N$}, up to logarithms (recall that $N=nT$
in the multi-task learning). The bounds (\ref{eq:bochk3}) and
(\ref{eq:bochk4}) achieve this rate in a more precise sense because
they are free of extra logarithmic factors.

Another remark concerns the possible range of $m$. It follows from the
discussion in Section \ref{sec: Schatten} that the ``dimension larger
than the sample size'' framework is not covered by Corollary
\ref{cor:ml2} since this corollary relies on the RI condition. In
contrast, the bounds of Corollary \ref{cor:ml3} make sense when the
dimension $m$ is larger than the sample size $n$ of each task; we only
need to have
$m\ll\exp(n)$ for Corollary \ref{cor:ml3} to be meaningful. 
Corollary \ref{cor:ml3} holds when the RI assumption is violated and
under a mild condition on the masks $X_i$. The price to pay is to
assume that the singular values of $A^*$ do not grow exponentially
fast. Also, the estimator of Corollary \ref{cor:ml3} corresponds to
$p<1$, so it is computationally hard.

\section{Minimax lower bounds}
\label{sec: lower bounds}
In this section, we derive lower bounds for the prediction error, which
show that the upper bounds that we have proved are optimal in a minimax
sense for two scenarios: (i) under the RI condition and (ii) for matrix
completion with collaborative sampling. We also provide a lower bound
for USR matrix completion. Under the RI condition with $\nu=1$, minimax
lower bounds for the Frobenius norm $\|\hat A -A^*\|_{S_2}$ on
``Schatten-0'' balls $\{A^*\in\R^{m\times T}\dvtx \operatorname{rank}(A^*)\leq
r\}$ are derived in Cand\`{e}s and Plan (\citeyear{cp09b}) with a technique
different from ours, which does not allow one to include further
boundedness constraints on $A^*$ in addition to that it has rank at
most $r$. Specifically, they prove their lower bound by passage to
Bayes risk with an unbounded support prior (Gaussian prior). Our lower
bounds are more general in the sense that they are obtained on smaller
sets, namely, the intersections of Schatten-0 and Schatten-$p$ balls.
This is similar in spirit to Rigollet and Tsybakov (\citeyear{rigtsy10}) establishing
minimax lower bounds on the intersection of $\ell_0$ and $\ell_1$ balls
for the vector sparsity scenario. In what follows, we denote by
$\inf_{\hat{A}}$ the infimum over all estimators based on
$(X_1,Y_1),\ldots,(X_N,Y_N)$, and for any $A\in\R^{m\times T}$ we denote
by $\PP_A$ the probability distribution of $(Y_1,\dots,Y_N)$ satisfying~(\ref{1}) with $A^*=A$.

\begin{theorem}[(Lower bound---Restricted Isometry)]\label{thm: lower
bound 1}
Let $\xi_1,\dots,\xi_N$ be i.i.d. $\NN(0,\sigma^2)$ random variables
for some $\sigma^2>0$. Let $M=\max(m,T)\geq8$, $r\geq1$,
$\min(T,m)\geq r$ and for $0<\alpha<1/8$ define
\begin{eqnarray*}
\psi_{M,N,r,\Delta}
&=&
\min\biggl(\frac{rM}{N}, \Delta^p\biggl(\frac{M}{N}\biggr)^{1-p/2},\Delta^2\biggr),
\\
C(\alpha)
&=&
\frac{\alpha(1-\delta_{r})^2}{2^{2/p}(1+\delta_{r})^2}\frac{\log2}{128}
\quad\mbox{and}\quad C(\alpha,\nu)=\frac{\alpha}{2^{2/p}(1+\delta_r)^2}\frac{\log2}{128}\nu^2.
\end{eqnarray*}
\begin{longlist}[(ii)]
\item[(i)] Assume that the sampling operator $\LL$ satisfies the
right-hand side
inequality in the RI~$(r,\nu)$-condition (\ref{eq:RIcondition}) for
some $\delta_{r}\in(0,1)$. Then for any $p\in(0,2]$, $\Delta>0$,
$0<\alpha<1/8$,
%
\begin{equation}\label{eq:lower2}
\qquad\inf_{\hat{A}}\mathop{\sup_{A^*\in\R^{m\times
T}:}}_{\operatorname{rank}(A^*)\leq r, \Vert A^*\Vert_{S_p}\leq\Delta\nu\sigma}
\PP_{A^*}\bigl( \|\hat A -A^*\|_{S_2}^2
> C(\alpha,\nu)\sigma^2\psi_{M,N,r,\Delta}\bigr)\geq\beta,
\end{equation}
where $\beta=\beta(M,\alpha)>0$ is such
that
$\beta(M,\alpha)\rightarrow1$ as $M\rightarrow\infty, \alpha\to0$.
\item[(ii)] Assume that the sampling operator $\LL$ satisfies the
RI~$(r,\nu)$-condition (\ref{eq:RIcondition}) for some
$\delta_{r}\in(0,1)$. Then for any $p\in(0,2]$, $\Delta>0$,
$0<\alpha<1/8$, with $\beta$ as in~(\ref{eq:lower2}),
%
\begin{equation}\label{eq:lower1}
\qquad\inf_{\hat{A}}\mathop{\sup_{A^*\in\R^{m\times
T}:}}_{\operatorname{rank}(A^*)\leq r, \Vert A^*\Vert_{S_p}\leq
\Delta\nu\sigma}\PP_{A^*}\bigl(\hat{d}_{2,N}(\hat{A},A^*)^2>
C(\alpha)\sigma^2\psi_{M,N,r,\Delta}\bigr)\geq\beta.
\end{equation}
\end{longlist}
\end{theorem}

\begin{ree}\label{re3} 
It is worth to note that $C(\alpha)$ and $\beta(M,\alpha)$ do not
depend on the constant $\nu$ of the RI condition.
\end{ree}

\begin{pf*}{Proof of Theorem \ref{thm: lower bound 1}}
Without loss of generality, we assume that
$M=m\geq T$.
For a constant $\gamma>0$ and an integer $s\in\{1,2,\ldots,r\}$, both to
be specified later, define
\[
\AAA_{s,\gamma}=\bigl\{ A=(a_{ij})\in\R^{m\times T}\dvtx a_{ij}\in\bigl\{0,
\gamma\nu/\sqrt{N}\bigr\} \mbox{ if } 1\leq j\leq s; a_{ij}=0
\mbox{ otherwise}\bigr\}.
\]
By construction, any element of $\AAA_{s,\gamma}$ as well as the
difference of any two elements of $\AAA_{s,\gamma}$ has rank at most
$s$. Due to the Varshamov--Gilbert bound [cf. Lemma~2.9 in
Tsybakov (\citeyear{tsyrr09})], there exists a subset
$\AAA^0_{s,\gamma}\subset\AAA_{s,\gamma}$ of cardinality
$\operatorname{Card}(\AAA_{s,\gamma}^0) \geq2^{sm/8}$ containing $A_0=0$ such
that for any two distinct elements $A_1$ and $A_2$ of
$\AAA_{s,\gamma}^0$,
%
\begin{equation}\label{eq: condition A}
\hat{d}_{2,N}(A_1,A_2)^2\geq\nu^{-2}(1-\delta_{r})^2\Vert
A_1-A_2\Vert_{S_2}^2\geq(1-\delta_{r})^2\frac{\gamma
^2}{8}\frac{sM}{N},
\end{equation}
where the first inequality follows from the left-hand side inequality in the
RI condition~(\ref{eq:RIcondition}) and is only used to prove
(\ref{eq:lower1}). We will prove (\ref{eq:lower1}); the proof of
(\ref{eq:lower2}) is analogous in view of the second inequality in
(\ref{eq: condition A}).

Then, for any $A_1\in\AAA_{s,\gamma}^0$, the Kullback--Leibler
divergence $K(\PP_{A_0},\PP_{A_1})$ between $\PP_{A_0}$ and $\PP_{A_1}$
satisfies
%
\begin{equation}\label{KLdiv}
K(\PP_{A_0},\PP_{A_1})=\frac{N}{2\sigma^2}\hat
{d}_{2,N}(A_0,A_1)^2\leq\frac{\gamma^2}{2\sigma^2}(1+\delta_{r})^2sM,
\end{equation}
where we used again the RI condition. We now apply Theorem 2.5 in
Tsybakov (\citeyear{tsyrr09}). Fix some $\alpha\in(0,1/8)$. Note that the condition
%
\begin{equation}\label{eq: condition C}
\frac{1}{\operatorname{Card}(\AAA_{s,\gamma}^0)-1}
\sum_{A\in\AAA_{s,\gamma}^0}K(\PP_A,\PP_{A_0})\leq\alpha\log
\bigl(\operatorname{Card}(\AAA_{s,\gamma}^0)-1\bigr)
\end{equation}
is satisfied for $ \gamma^2 \leq\alpha\sigma^2(\log
2)/(4(1+\delta_{r})^2)$. Define
\[
r_{\Delta}=\operatorname{arg\,min}\biggl\{l\in\N\dvtx \Delta^p\leq l\biggl(\frac{M}{N}\biggr)^{p/2}\biggr\},
\]
and consider separately the following three cases.

\textit{The case $r_{\Delta}=1$.} In this case,
$\psi_{M,N,r,\Delta}=\Delta^2$ for any $r\geq1$, and $\Delta^2N/M\leq 1$. Set
\[
s_1=1\quad\mbox{and}\quad\gamma_1=\biggl(\frac{\alpha}{(1+\delta_{r})^2}\frac
{\log
2}{4}\sigma^2\Delta^2\frac{N}{M}\biggr)^{1/2}.
\]
Then $\Vert A\Vert_{S_p} \leq\Vert A\Vert_{S_2}
\leq\sqrt{M/N}\nu\gamma\leq\Delta\nu\sigma$ for all
$A\in\AAA_{1,\gamma_1}$, i.e., $\AAA_{1,\gamma_1}$ is contained in the
set
\[
\{A\in\R^{m\times T}\dvtx \operatorname{rank}(A)\leq r,\Vert A\Vert
_{S_p}\leq\Delta\nu\sigma\}.
\]
Now, inequality (\ref{eq: condition A}) shows that $
\hat{d}_{2,N}(A_1,A_2)^2 \geq4C(\alpha)\sigma^2\Delta^2$ for any two
distinct elements $A_1,A_2\in\AAA_{1,\gamma_1}^0$, while
$\Delta^2N/M\leq1$ implies that $\gamma_1^2\leq\break\alpha\sigma ^2(\log
2)/(4(1+\delta_{r})^2)$. Hence, condition (\ref{eq: condition C}) is
satisfied.

\textit{The case $2\!\le\!r_{\Delta}\!\leq\!r$.} In this
case, the rate $\psi_{M,N,r,\Delta}$ is equal to
$\Delta^p(M/N)^{1-p/2}$. We consider the set
$\AAA_{r_{\Delta},\gamma_2}^0$ with some $\gamma_2$ to be specified
below. For $A\in\AAA_{r_{\Delta},\gamma_2}^0$, we have $\Vert
A\Vert_{S_p}^2\leq r_{\Delta}^{(2-p)/p}\Vert
A\Vert_{S_2}^2 \leq r_{\Delta}^{2/p}\gamma_2^2\nu^2 M/N$. Since
also $r_{\Delta}\leq2 \Delta^p(N/M)^{p/2}$ when $r_{\Delta}\ge2$, it
follows that $\Vert A\Vert_{S_p}\leq\Delta\nu\sigma$ whenever
%
\begin{equation}\label{eq: hi 2}
\gamma_2\leq2^{-1/p}\sigma.
\end{equation}
%
Now define
\[
s_2=r_{\Delta}\quad\mbox{and}\quad
\gamma_2=2^{-1/p}\biggl(\frac{\alpha}{(1+\delta_{r})^2}\frac{\log
2}{4}\sigma^2\biggr)^{1/2}.
\]
Then (\ref{eq: condition C}) is satisfied and $\gamma_{2}$ fulfills
also the constraint (\ref{eq: hi 2}), since $\alpha<1/8$, $(\log2)/4<
1$. Thus, $\AAA_{r_{\Delta},\gamma_2}^0$ is a subset of matrices
$A\in\R^{m\times T}$ with $\operatorname{rank}(A)\leq r$ and $\Vert
A\Vert_{S_p}\leq\Delta\nu\sigma$. Finally, (\ref{eq: condition A})
implies that
\begin{eqnarray*}
\hat{d}_{2,N}(A_1,A_2)^2 &\geq&
(1-\delta_r)^2\frac{\gamma_2^2}{8}\frac{r_{\Delta}M}{N} \geq
(1-\delta_r)^2\frac{\gamma_2^2}{8}\Delta^p\biggl(\frac{M}{N}\biggr)^{1-p/2}
\\
 &=& 4C(\alpha)\sigma^2\Delta^p\biggl(\frac{M}{N}\biggr)^{1-p/2}
\end{eqnarray*}
for any two distinct elements $A_1,A_2$ of
$\AAA_{r_{\Delta},\gamma_2}$.

\textit{The case $r_{\Delta} > r$.} In this case,
$\psi_{M,N,\Delta,r}=rM/N$. The conditions required in Theorem 2.5 of
Tsybakov (\citeyear{tsyrr09}) follow immediately as above, this time with the set of
matrices $\AAA_{r,\gamma_3}^0$, where
$\gamma_3^2=\alpha\sigma^2(\log2)/(4(1+\delta_r)^2)$.
\end{pf*}

\begin{ree}
Theorem~\ref{thm: lower bound 1}
implies that the rates of convergence in Theorem~\ref{th:00} are
optimal in a minimax sense on Schatten-$p$ balls $\{A^*\in\R^{m\times
T}\dvtx\break\|A^*\|_{S_p}\leq\Delta\}$ under the RI condition and natural
assumptions on $m,T$ and $N$. Indeed, using Theorem~\ref{thm: lower
bound 1} with no restriction on the rank [i.e., when $r=\min(m,T)$],
and putting for simplicity $\Delta=1$, we find that the rate in the
lower bound is of the order $\min(\min(m,T)M/N, (M/N)^{1-p/2}, 1)$. For
$m=T(=M)$ and $m^3>N>m$ this minimum equals $(M/N)^{1-p/2}$, which
coincides with the upper bound of Theorem~\ref{th:00}.
\end{ree}

The lower bound for the prediction error (\ref{eq:lower1}) in the above
theorem does not apply to matrix completion with $N<mT$ since then the
Restricted Isometry condition cannot be satisfied, as discussed in
Section \ref{sec: Schatten}. However, for the bound (\ref{eq:lower2})
we only need the right-hand side inequality in the RI condition. For
example, the latter is trivially satisfied for CS matrix completion
with $\nu=\sqrt{N}$ and $\delta_r=0$. This yields the following
corollary.

\begin{corollary}[(Lower bound---CS matrix completion)]
\label{cor:lower_collabor} Let $\xi_1,\dots,\xi_N$ be i.i.d.
$\NN(0,\sigma^2)$ random variables for some $\sigma^2>0$. Let
$M=\max(m,T)\geq8$, $r\geq1$, $\min(T,m)\geq r$, and consider the
problem of CS matrix completion.
Then for any
$p\in(0,2]$, $\Delta>0$, $0<\alpha<1/8$,
\[
\inf_{\hat{A}}
\mathop{\sup_{A^*\in\R^{m\times
T}:}}_{\operatorname{rank}(A^*)\leq r, \Vert A^*\Vert_{S_p}\leq
\Delta\sqrt{N}\sigma}\PP_{A^*}\biggl( \frac{1}{N}\|\hat A
-A^*\|_{S_2}^2
> C'(\alpha)\sigma^2\psi_{M,N,r,\Delta}\biggr)\geq\beta,
\]
where $C'(\alpha)=\alpha(\log2)2^{-2/p}/128$ and
$\beta=\beta(M,\alpha)$,
$\psi_{M,N,r,\Delta}$ are as in Theorem~\ref{thm: lower bound 1}.
\end{corollary}

The model of uniform sampling without replacement considered in
Cand\`{e}s and Recht (\citeyear{recht09}) is a particular case of CS matrix
completion. In the noisy case, Keshavan, Montanari and Oh (\citeyear{kmoff09}) obtain upper
bounds under such a sampling scheme with the rate $rM/N$, up to
logarithmic factors. The lower bound of
Corollary~\ref{cor:lower_collabor} is of the same order when
$\Delta=\infty$, that is, for the class of matrices of rank smaller
than $r$. However, Keshavan,
Montanari and Oh (\citeyear{kmoff09}) obtained their bounds on some
subclasses of this class characterized by additional strong
restrictions.

It is useful to note that for bounds of the type (\ref{eq:lower2}) it
is enough to have a condition on $\LL$ in expectation, as specified in
the next theorem.

\begin{theorem}
\label{thm:lowerbound_USR}
Let $\xi_1,\dots,\xi_N$ be i.i.d.
$\NN(0,\sigma^2)$ random variables for some $\sigma^2>0$. Let
$M=\max(m,T)\geq8$, $r\geq1$, $\min(T,m)\geq r$, and assume that
$X_1,\dots,X_N$ are random matrices independent of $\xi_1,\dots,\xi_N$,
and the sampling operator satisfies $\nu^2 \E|{\mathcal L}(A)|_2^2
\le
\|A\|_{S_2}^2$ for some $\nu>0$ and all $A\in\R^{m\times T}$ such that
${\rm rank} (A)\le r$. Then for any $p\in(0,2]$, $\Delta>0$,
$0<\alpha<1/8$,
\[
\inf_{\hat{A}}\mathop{\sup_{A^*\in\R^{m\times
T}:}}_{\operatorname{rank}(A^*)\leq r, \Vert A^*\Vert_{S_p}\leq
\Delta\nu\sigma}\PP_{A^*}\biggl(\frac{1}{\nu^2} \|\hat A
-A^*\|_{S_2}^2
> C'(\alpha)\sigma^2\psi_{M,N,r,\Delta}\biggr)\geq\beta,
\]
where $C'(\alpha)=\alpha(\log2)2^{-2/p}/128$ and
$\beta=\beta(M,\alpha)$,
$\psi_{M,N,r,\Delta}$ are as in Theorem~\ref{thm: lower bound 1}.
\end{theorem}

\begin{pf}
We proceed as in Theorem~\ref{thm: lower bound
1}, with the only difference in the bound on the Kullback--Leibler
divergence. Indeed, under our asumptions, instead of (\ref{KLdiv}) we
have
%
\begin{equation}\label{KLdiv1}
\qquad K(\PP_{A_0},\PP_{A_1})=\frac{N}{2\sigma^2}\E(\hat
{d}_{2,N}(A_0,A_1)^2)\le\frac{N}{2\nu^2}\|A_0
-A_1\|_{S_2}^2\le\frac{\gamma^2sM}{2\sigma^2}.
\end{equation}
\upqed\end{pf}

Theorem~\ref{thm:lowerbound_USR} applies to USR matrix completion with
$\nu=\sqrt{mT}$. Indeed, in that case $mT \E|{\mathcal L}(A)|_2^2 =
\|A\|_{S_2}^2$. In particular, Theorem~\ref{thm:lowerbound_USR} with
$\Delta=\infty$ shows that on the class of matrices of rank smaller
than $r$ the lower bound of estimation in the squared Frobenius norm
for USR matrix completion is of the order $rM/N$.

The next theorem gives a lower bound for the prediction error under
collaborative sampling without the RI condition. Instead, we only
impose a rather natural condition that the observed noisy entries are
sufficiently well dispersed, that is, there exist $r$ rows or $r$
columns with more that $\kappa Mr$ observations for some fixed
$\kappa\in(0,1]$. We state the result with an additional constraint on
the Frobenius norm of $A^*$, in order to fit the corresponding upper
bound (cf. Remark 2 in Section \ref{sec:examp}).

\begin{theorem}[(Lower bound---CS matrix completion)]\label{thm: lower
bound 2}
Let $\xi_1,\dots,\xi_N$ be i.i.d. $\NN(0,\sigma^2)$ random variables
for some $\sigma^2>0$ and assume that the masks
$X_1=e_{i_1}(m)e_{j_1}'(T),\dots, X_N=e_{i_N}(m)e_{j_N}'(T)$ are
pairwise different. Let $\min(T,\break m)\geq r$ and $\kappa Mr\geq8$ for
some fixed $\kappa\in(0,1]$, where $M=\max(m,T)$. Assume furthermore
that the following dispersion condition holds: there exist numbers
$1\leq k_1<\cdots<k_r\leq T$ or $1\leq k_1'<\cdots<k_r'\leq m$ such that
either the set
$\{(i_1,j_1),\dots,(i_N,j_N)\}\,\cap\,\{(i,k_1),\ldots,(i,k_r)\dvtx i=1,\ldots,m\}$
or the set $\{(i_1,j_1),\ldots,(i_N,j_N)\}\,\cap\,\{(k_1',j),\ldots,(k_r',j)\dvtx
j=1,\ldots,T\} $ has cardinality at least $\kappa Mr+1$. Define $
\CC_{\delta,r} = \{A\in\R^{m\times T}\dvtx \operatorname{rank}(A)\leq
r\mbox{ and }\Vert A\Vert_{S_2}\leq\delta\}. $ Then for any
$0<\alpha<1/8$ and $ \delta^2\geq\alpha\sigma^2(\log2)(\kappa M
r+1)/4$,
\[
\inf_{\hat{A}}\sup_{A^*\in\CC_{\delta,r}}\PP_{A^*}\biggl(\hat
{d}_{2,N}(\hat{A},A^*)^2> C'(\alpha)\frac{\sigma^2\kappa
rM}{N}\biggr)\geq\beta(\kappa M,\alpha)>0,
\]
with a function $\beta\rightarrow1$ as $\kappa M\rightarrow\infty,
\alpha\rightarrow0$ and $C'(\alpha)=\alpha(\log2)/128$.
\end{theorem}

\begin{pf}
We proceed as in the case $\Delta=\infty$,
$p=2$, $\nu=\sqrt{N}$
of Theorem~\ref{thm: lower bound 1} taking a different set ${\mathcal
A}^0$ instead ${\mathcal A}^0_{s,\gamma}$. Let, for definiteness, the
dispersion condition be satisfied with the set of indices ${\mathcal
K}= \{(i_1,j_1),\dots,(i_N,j_N)\}\cap\{(i,k_1),\dots,(i,k_r)\dvtx i=1,\ldots,m\}$. Then there exists a subset ${\mathcal
K}'$ of ${\mathcal
K}$ with cardinality $\operatorname{Card}({\mathcal
K}')=\lceil\kappa Mr \rceil$. We define
\[
\AAA=\bigl\{ A=(a_{ij})\in\R^{m\times T}\dvtx a_{ij}\in\{0, \gamma\} \mbox{ if }
(i,j)\in{\mathcal K}'; a_{ij}=0 \mbox{ otherwise}\bigr\}.
\]
Any element of $\AAA$ as well as the difference of any two elements of
$\AAA$ has rank at most $r$, and $\Vert A\Vert_{S_2}^2\le
\gamma^2 \lceil\kappa Mr \rceil$, $\forall A\in\AAA$. So, $\AAA \subset
\CC_{\delta,r}$ if $\gamma^2 ( \kappa Mr +1)\le\delta^2$. As in
Theorem~\ref{thm: lower bound 1}, the Varshamov--Gilbert bound implies
that there exists a~subset $\AAA^0\subset\AAA$ of cardinality
$\operatorname{Card}(\AAA^0) \geq2^{\lceil\kappa Mr \rceil/8}$ containing
$A_0=0$, that for any two distinct elements $A_1$ and $A_2$ of
$\AAA^0$,
\[
\hat{d}_{2,N}(A_1,A_2)^2= N^{-1}\Vert
A_1-A_2\Vert_{S_2}^2\geq\frac{\gamma^2}{8}\frac{\lceil\kappa Mr
\rceil}{N}.
\]
Instead of the bound (\ref{KLdiv}), we have now the inequality
$K(\PP_{A_0},\PP_{A_1}) \leq\frac{\gamma^2}{2\sigma^2} \lceil \kappa Mr
\rceil$ for any $A_1\in\AAA^0$. Finally, we choose $\gamma^2=
\alpha\sigma^2 (\log2)/4$. With these modifications, the rest of the
proof is the same as that of Theorem~\ref{thm: lower bound 1} in the
case $r_{\Delta}>r$.
\end{pf}

\section{Control of the stochastic term}
\label{sec: stochastic term}

We consider two approaches for bounding the stochastic term $
N^{-1}\sum_{i=1}^N\xi_i\operatorname{tr}((\hat{A}-A^*)'X_i)
$ on the right-hand side of the basic inequality (\ref{eq: BI}).
The first one used for $p=1$ consists in application of the trace
duality
%
\begin{equation}\label{td}
\Biggl\vert
\frac{1}{N}\sum_{i=1}^N\xi_i\trace\bigl((\hat{A}-A^*)'X_i\bigr)\Biggr\vert
\leq\Vert
\hat{A}-A^*\Vert_{S_1}\Vert\mathbf{M}\Vert
_{S_{\infty}}
\end{equation}
with $\mathbf{M}=N^{-1}\sum_{i=1}^N\xi_i X_i$ and then of suitable
exponential bounds for the spectral norm of $\mathbf{M}$ under different
conditions on $X_i$, $i=1,\ldots,N$. The second approach used to treat the
case $0<p<1$ (nonconvex penalties) (cf. Section~\ref{subsec: EP}) is
based on refined empirical process techniques. Proofs of the results of
this section are deferred to Section~\ref{sec:proof_lem}.

\subsection{Tail bounds for the spectral norm of random matrices}
\label{subsec: tail bounds}
We say that the random variables $\xi_i$, $i=1,\dots,N$, satisfy the
\textit{Bernstein condition} if
%
\begin{equation}\label{eq: Bernstein}
\max_{1\le i\le
N}\E\vert\xi_i\vert^l\leq\frac{1}{2}l!\sigma^2H^{l-2},\qquad l=
2,3,\ldots,
\end{equation}
with some finite constants $\sigma$ and $H$, and we say that they
satisfy the \textit{light tail condition} if
%
\begin{equation}\label{eq: LT}
\max_{1\leq i\leq N}\E(\exp(\xi_i^2/\sigma^2))\leq\exp(1)
\end{equation}
for some positive constant $\sigma^2$.

\begin{lemma}\label{lemma:bernstein}
Let the i.i.d. zero-mean random variables $\xi_i$ satisfy the Bernstein
condition (\ref{eq: Bernstein}). Let also either
%
\begin{equation}\label{eq: con 2a}
\max_{1\leq j\leq m}\frac{1}{N}\sum_{i=1}^N\bigl\vert
X_{i(j,\cdot)}\bigr\vert_2^2\leq S_{\mathrm{row}}^2
\end{equation}
and
\begin{equation}\label{eq: con 2b}
\max_{1\le j\leq m, 1\le i\leq N}\bigl\vert
X_{i(j,\cdot)}\bigr\vert_2\leq H_{\mathrm{row}}
\end{equation}
or the conditions
%
\begin{equation}\label{eq: con 1a}
\max_{1\leq k\leq T}\frac{1}{N}\sum_{i=1}^N\bigl\vert X_{i(\cdot
,k)}\bigr\vert_2^2\leq S_{\mathrm{col}}^2
\end{equation}
and
\begin{equation}\label{eq: con 1b}
\max_{1\le k\leq T, 1\le i\leq N}\bigl\vert
X_{i(\cdot,k)}\bigr\vert_2\leq H_{\mathrm{col}}
\end{equation}
hold true with some constants $S_{\mathrm{row}},H_{\mathrm{row}}, S_{\mathrm{col}},H_{\mathrm{col}}$. Let
$D>1$. Then, respectively, with probability at least $1-2/m^{D-1}$ or
at least $1-2/T^{D-1}$ we have
%
\begin{equation}\label{eq:stochterm1}
\Vert\mathbf{M} \Vert_{S_{\infty}}\leq\tau,
\end{equation}
where $\tau=\tau_{\mathrm{row}}=C_{\mathrm{row}}\sqrt{m(\log m)/N}$ if (\ref{eq: con 2a})
and (\ref{eq: con 2b}) are satisfied or $\tau=\tau_{\mathrm{col}}=
C_{\mathrm{col}}\sqrt{T(\log T)/N}\}$ if (\ref{eq: con 1a}) and (\ref{eq: con
1b}) hold. Here
\begin{eqnarray*}
C_{\mathrm{row}}&=& \Biggl(\sqrt{2 D\sigma^2S_{\mathrm{row}}^2} + 2DH_{\mathrm{row}}H\sqrt{\frac{\log
m}{N}}\Biggr),
\\
 C_{\mathrm{col}}&=& \Biggl(\sqrt{2 D\sigma^2S_{\mathrm{col}}^2} +
2DH_{\mathrm{col}}H\sqrt{\frac{\log T}{N}}\Biggr).
\end{eqnarray*}
\end{lemma}


\begin{lemma}\label{lemma:stochterm_matr_compl_1}
Let $\xi_1,\ldots,\xi_N$ be i.i.d. $\NN(0,\sigma^2)$ random variables.
Then, for any $D\ge2$,
%
\begin{equation}\label{eq:stochterm00}
\Vert\mathbf{M} \Vert_{S_{\infty}}\leq4\sqrt{2D}\sigma
\phi_{\max}(1)\sqrt{\frac{m+T}{N}}=:\tau_1
\end{equation}
with probability at least $1-2\exp\{-(D-\log
5)(m+T)\}$, where $\phi_{\max}(1)$ is the maximal rank 1 eigenvalue of
the sampling operator $\LL$.
\end{lemma}


If $m$ and $T$ have the same order of magnitude, the bound of
Lemma~\ref{lemma:stochterm_matr_compl_1} is better, since it does not
contain extra logarithmic factors. On the other hand, if $m$ and $T$
differ dramatically, for example, $m\gg T$, then
Lemma~\ref{lemma:bernstein} can provide a significant improvement.
Indeed, the ``column'' version of Lemma~\ref{lemma:bernstein} guarantees
the rate $\tau\sim\sqrt{T\log T}/\sqrt{N}$ which in this case is much
smaller than $\sqrt{m/N}$. In all the cases, the concentration rate in
Lemma~\ref{lemma:stochterm_matr_compl_1} is exponential and thus faster
than in Lemma~\ref{lemma:bernstein}.


The next lemma treats the stochastic term for USR matrix completion.

\begin{lemma}[(USR matrix completion)]\label{lemma:stochterm_matr_compl}
\textup{(i)} Let the i.i.d. zero-mean random variables $\xi_i$ satisfy the
Bernstein condition (\ref{eq: Bernstein}). Consider the USR matrix
completion problem and assume that $mT(m+T)>N$. Then, for any $D\ge2$,
%
\begin{equation}\label{eq:stochterm0}
\Vert\mathbf{M}\Vert_{S_{\infty}}\leq\bigl(4\sigma\sqrt
{10D} +
8HD\bigr)\frac{m+T}{N}=:\tau_2
\end{equation}
with probability at least $1-4\exp\{-(2-\log5)(m+T)\}$.

\textup{(ii)} Assume that the i.i.d. zero-mean random variables $\xi_i$ satisfy
the light tail condition (\ref{eq: LT}) for some $\sigma^2>0$. Then for
any $B>0$,
%
\begin{equation}\label{eq:nemirov}
\Vert
\mathbf{M}\Vert_{S_{\infty}}\leq\sqrt{B}\frac{\sigma\log
(\max(m+1,T+1))}{\sqrt{N}}=:
\tau_3
\end{equation}
with probability at least $1-(1/C){\max(m+1,T+1)}^{-CB}$ for some
constant $C>0$ which does not depend on $m,T$ and $N$.
\end{lemma}

%

The proof of part (i) is based on a refinement of a technique in
Vershynin (\citeyear{vershf07}), whereas that of part (ii) follows immediately from
the large deviations inequality of Nemirovski~(\citeyear{nem04}). For example, if
$\xi_i\sim\NN(0,\sigma^2)$, in which case both results apply, the bound
(ii) is tighter than (i) for sample sizes $N\ll(m+T)^2$ which is the
most interesting case for matrix completion.

Much tighter bounds are available when the $X_i$ are constrained to be
pairwise different. Besides it is noteworthy that the rates in
(\ref{eq:stochterm010}) and (\ref{eq:stochterm011}) below are different
for Gaussian and Bernstein errors.

\begin{lemma}[(Collaborative sampling)]\label{lemma:stochterm_matr_compl_2}
Consider the problem of CS matrix completion.
\begin{longlist}[(iii)]
\item[(i)] Let $\xi_1,\ldots,\xi_N$ be i.i.d. $\NN(0,\sigma^2)$ random variables.
Then, for any\break $D\ge2$,
%
\begin{equation}\label{eq:stochterm010}
\Vert\mathbf{M}\Vert_{S_{\infty}}
\leq8\sigma\sqrt{D}\frac{\sqrt{m+T}}{N}=:\tau_4
\end{equation}
with probability at least $1-2\exp\{-(D-\log5)(m+T)\}$.
\item[(ii)] Let $\xi_1,\ldots,\xi_N$ be i.i.d. zero-mean random variables
satisfying the Bernstein condition (\ref{eq: Bernstein}). Then, for any
$D\ge2$ and
%
\begin{equation}\label{eq:stochterm011}
\Vert\mathbf{M}\Vert_{S_{\infty}}
\leq\frac{4\sigma\sqrt{2D(m+T)} + 8HD(m+T)}{N}=:\tau_5
\end{equation}
with probability at least $1-2\exp\{-(D-\log5)(m+T)\}$.
\item[(iii)] Let $\xi_1,\ldots,\xi_N$ be i.i.d. $\NN(0,\sigma^2)$ random
variables. Then for any $A>1$,
\[
\Vert\mathbf{M}\Vert_{S_{\infty}}
\leq\frac{\sigma\sqrt{2A\log(m + T)}}{N}\max\Biggl\{ \Biggl\Vert
\sum_{i=1}^NX_i'X_i\Biggr\Vert_{S_{\infty}}^{1/2}, \Biggl\Vert
\sum_{i=1}^NX_iX_i'\Biggr\Vert_{S_{\infty}}^{1/2} \Biggr\} =:\tau_6
\]
with probability at least $1-2(m+T)^{1-A}$.
\end{longlist}
\end{lemma}


Since the masks $X_i$ are distinct, the maximum appearing in (iii) is
bounded by $\sqrt{\max(m,T)}$; in case it is attained, the bound
(\ref{eq:stochterm010}) is slightly stronger since it is free from the
logarithmic factor. For $N\ll mT$ the tightness of the bound in (iii)
depends strongly on the geometry of the $X_i$'s and the
maximum can be significantly smaller than $\sqrt{\max(m,T)}$. 
Note also that the concentration in (\ref{eq:stochterm010}) is
exponential, while it is only polynomial in (iii).

\subsection{Concentration bounds for the stochastic term under
nonconvex penalties} \label{subsec: EP}
The last bound in this section applies in the case $0<p<1$. It is given
in the following lemma.

\begin{lemma}\label{nonconv}
Let $\xi_1,\ldots,\xi_N$ be i.i.d. $\NN(0,\sigma^2)$ random variables,
$0<p<1$ and $M=\max(m,T)$. Assume that the sampling operator $\LL$ is
uniformly bounded; cf. (\ref{eq: UB}). Set
$c_{\kappa}=(2\kappa-1)(2\kappa)\kappa^{-1/(2\kappa-1)}$ where
$\kappa=(2-p)/(2-2p)$. Then for any fixed $\delta>0$, $\vartheta\geq
C^2$ and $ \tau_7= c_{\kappa}(\vartheta/p)^{1-p/2}(M/N)^{1-p/2} $ we
have
%
\begin{equation}\label{eq:nonconv}
\quad\quad\Biggl\vert\frac{1}{N}\sum_{i=1}^N\xi_i\operatorname{tr}\bigl(X_i'(\hat
{A}-A^*)\bigr)\Biggr\vert\leq\frac{\delta}{2}\hat{d}_{2,N}(\hat{A},A^*)^2
+ \tau_7\delta^{p-1}\Vert\hat{A}-A^*\Vert_{S_p}^{p}\hspace*{-10pt}
\end{equation}
with probability at least $1- C\exp(-\vartheta M/C^2)$ for some
constant $C=C(p,c_0,\break\sigma^2)>0$ which is independent of $M$ and $N$
and satisfies $ \sup_{0<p\leq q}C(p,c_0,\sigma)<\infty$ for all $q<1$.
\end{lemma}

Note at this point that we cannot rely the proof of Lemma \ref{nonconv}
directly on the trace duality and norm interpolation (cf. Lemma
\ref{lemma7a}), that is, on the inequalities
\begin{eqnarray}\label{chain1}
\quad\quad \Biggl\vert
\frac{1}{N}\sum_{i=1}^N\xi_i\operatorname{tr}\bigl(X_i'(\hat
{A}-A^*)\bigr)\Biggr\vert&\leq&\Vert
\hat{A}-A^*\Vert_{S_1}\Vert\mathbf{M}\Vert_{S_{\infty}}\nonumber
\\[-8pt]\\[-8pt]
&\leq&\Vert
\hat{A}-A^*\Vert_{S_2}^{1-p/(2-p)}\Vert
\hat{A}-A^*\Vert_{S_p}^{p/(2-p)}\Vert\mathbf{M}\Vert_{S_{\infty}}\nonumber.\hspace*{-10pt}
\end{eqnarray}
Indeed, one may think that we could have bounded here the
$S_\infty$-norm of $\mathbf{M}$ in the same way as in Section \ref{subsec:
tail bounds}, and then the proof would be complete after suitable
decoupling if we were able to bound from above
$\Vert\hat{A}-A^*\Vert_{S_2}^2$ by
$\hat{d}_{2,N}(\hat{A},A^*)^2$ times a constant factor. However, this
is not possible. Even the Restricted Isometry condition cannot help
here because $\hat{A}-A^*$ is not necessarily of small rank.
Nevertheless, we will show that by other techniques it is possible to
derive an inequality similar to (\ref{chain1}) with
$\hat{d}_{2,N}(\hat{A},A^*)$ instead of
$\Vert\hat{A}-A^*\Vert_{S_2}$. Further details are given in
Sections \ref{sec:proof_lem} and \ref{sec: embeddings}.

\section{\texorpdfstring{Proof of Theorem \protect\ref{main}}{Proof of Theorem 2}}
\label{sec: prooft1}
\subsection*{Preliminaries}
We first give two lemmas on matrix
decomposition needed in our
proof, which are essentially provided by Recht, Fazel and Parrilo
(\citeyear{rfp07}) [subsequently, RFP(10) for short].

\begin{lemma}\label{lemma: 2.3}
Let $A$ and $B$ be matrices of the same dimension. If $AB'=0$, $A'B=0$,
then
\[
\Vert A+B\Vert_{S_p}^p=\Vert
A\Vert_{S_p}^p+\Vert B\Vert_{S_p}^p\qquad \forall p>0.
\]
\end{lemma}

\begin{pf}
For $p=1$ the result is Lemma 2.3 in RFP(10).
The argument obviously extends to any $p>0$ since RFP(10) show that the
singular values of $A+B$ are equal to the union (with repetition) of
the singular values of $A$ and $B$.
\end{pf}

\begin{lemma}\label{lemma: 3.4}
Let $A\in\R^{m\times T}$ with $\operatorname{rank}(A)=r$ and singular value
decomposition $A=U\Lambda V'$. Let ${B}\in\R^{m\times T}$ be arbitrary.
Then there exists a decompositon $B=B_1+B_2$ with the following
properties:
\begin{longlist}[(iii)]
\item[(i)] $\operatorname{rank}(B_1)\leq2\operatorname{rank}(A)=2r$,
\item[(ii)] $A B_2'=0$, $A'B_2=0$,
\item[(iii)] $\operatorname{tr}(B_1'B_2) = 0$,
\item[(iv)] $B_1$ and $B_2$ are of the form
%
\begin{eqnarray}
&&B_1=U\pmatrix{
\tilde{B}_{11} &
\tilde{B}_{12}\cr
\tilde{B}_{21} & 0
}
V'\quad\mbox{and}\quad B_2=U\pmatrix{
0 & 0\cr
0 & \tilde{B}_{22}
}V'\nonumber
\\
&&\eqntext{\mbox{with }\tilde{B}_{11}\in\R^{r\times r}.}
\end{eqnarray}
\end{longlist}
\end{lemma}

The points (i)--(iii) are the statement of Lemma 3.4 in
RFP(08), the representation (iv) is provided in its proof.

\begin{pf*}{Proof of Theorem \ref{main}}
First note that there exists a decomposition
$\hat{A}=\hat{A}^{(1)}+\hat{A}^{(2)}$ with the following properties:
\begin{longlist}[(iii)]
\item[(i)] $\operatorname{rank}(\hat{A}^{(1)}-A^*)\leq
2\operatorname{rank}(A^*)=2r$,
\item[(ii)] $A^*(\hat{A}^{(2)})'=0$, $(A^*)'\hat{A}^{(2)}=0$,
\item[(iii)] $\operatorname{tr}((\hat{A}^{(1)}-A^*)'\hat
{A}^{(2)}) =
0$.
\end{longlist}
This follows from Lemma \ref{lemma: 3.4} with $A=A^*$ and
$B=\hat{A}-A^*$. In the notation of Lemma \ref{lemma: 3.4}, we have
$B_1=\hat{A}^{(1)}-A^*$ and
$B_2=\hat{A}^{(2)}$.

From the basic inequalities (\ref{eq: BI}) and (\ref{stoch}) with
$\delta=1/2$, we find
%
\begin{eqnarray}\label{eq: starting0}
&&\bigl(1-I_{\{0<p<1\}}/2\bigr)\hat{d}_{2,N}(\hat{A},A^*)^2\nonumber
\\[-8pt]\\[-8pt]
&&\qquad\leq2^{2-p}\tau
\Vert\hat{A}-A^*\Vert_{S_p}^p+4\tau(\Vert
A^*\Vert_{S_p}^p-\Vert\hat{A}\Vert_{S_p}^p).\nonumber
\end{eqnarray}
In particular, for the case $p=1$,
%
\begin{equation}\label{eq: starting}
\hat{d}_{2,N}(\hat{A},A^*)^2\leq2\tau\Vert\hat
{A}-A^*\Vert_{S_p}^p+4\tau(\Vert
A^*\Vert_{S_p}^p-\Vert\hat{A}\Vert_{S_p}^p).
\end{equation}
For brevity, we will conduct the proof with the numerical constants
given in (\ref{eq: starting}), that is, with those for $p=1$. The proof
for general $p$ differs only in the values of the constants, but their
expressions become cumbersome.

Using (\ref{triang}), we get
%
\begin{eqnarray}\label{eq: starting1}
&&\hat{d}_{2,N}(\hat{A},A^*)^2\nonumber
\\[-8pt]\\[-8pt]
&&\qquad\leq2\tau\bigl\Vert\hat
{A}^{(1)}-A^*\bigr\Vert_{S_p}^p+4\tau\Vert
A^*\Vert_{S_p}^p+2\tau\bigl\Vert\hat{A}^{(2)}\bigr\Vert_{S_p}^p-
4\tau\Vert\hat{A}\Vert_{S_p}^p.\nonumber
\end{eqnarray}
By (\ref{triang}) again and by Lemma \ref{lemma: 2.3},
\begin{eqnarray*}
\Vert\hat{A}\Vert_{S_p}^p
&\geq&
\bigl\Vert
A^*+\hat{A}^{(2)}\bigr\Vert_{S_p}^p-\bigl\Vert\hat{A}^{(1)}-A^*\bigr\Vert_{S_p}^p
\\
&=&
\Vert A^*\Vert_{S_p}^p+\bigl\Vert\hat{A}^{(2)}\bigr\Vert_{S_p}^p-\bigl\Vert\hat{A}^{(1)}-A^*\bigr\Vert_{S_p}^p,
\end{eqnarray*}
since $(A^*)'\hat{A}^{(2)}=0$ and $A^*(\hat{A}^{(2)})'=0$ by
construction. Together with (\ref{eq: starting1}) this yields
%
\begin{equation}
\qquad\hat{d}_{2,N}(\hat{A},A^*)^2\leq2\tau\bigl\Vert
\hat{A}^{(1)}-A^*\bigr\Vert_{S_p}^p-2\tau\bigl\Vert\hat
{A}^{(2)}\bigr\Vert_{S_p}^p+4\tau\bigl\Vert\hat
{A}^{(1)}-A^*\bigr\Vert_{S_p}^p,\hspace*{-5pt}
\end{equation}
from which one may deduce
%
\begin{equation}\label{eq: ub}
\hat{d}_{2,N}(\hat{A}, A^*)^2\leq6\tau\bigl\Vert
\hat{A}^{(1)}-A^*\bigr\Vert_{S_p}^p
\end{equation}
and
\begin{equation}\label{eq: 1--2}
\bigl\Vert\hat{A}^{(2)}\bigr\Vert_{S_p}^p\leq3\bigl\Vert
\hat{A}^{(1)}-A^*\bigr\Vert_{S_p}^p.
\end{equation}
Consider now the following decomposition of the matrix $\hat{A}^{(2)}$.
First, recall that $\hat{A}^{(2)}$ is of the form
\[
\hat{A}^{(2)}=U\pmatrix{ 0 & 0\cr0 & \tilde{B}_{22}
}V'.
\]
Write $\tilde{B}_{22}=W_1\Lambda(\tilde{B}_{22}) W_2'$ with diagonal
matrix $\Lambda(\tilde{B}_{22})$ of dimension $r'$ and
$W_1'W_1=W_2'W_2=I_{r'\times r'}$ for some $r'\le\min(m,T)$. In the
next step, $W_1$ and $W_2$ are complemented to orthogonal matrices
$\bar{W}_1$ and $\bar{W}_2$ of dimension $\min(m,T)\times\min(m,T)$.
For instance, set
\[
\bar{W}_2'=\pmatrix{
& & 0\cr
* & &\cr
& & W_2'
}\in\R^{\min(m,T)\times\min(m,T)},
\]
where $*$ complements the columns of the matrix ${{0}\choose{W_2'}}$ to an
orthonormal basis in $\R^{m\times T}$, and proceed analogously with
$W_1$. In particular,
$\bar{W}_1'\bar{W}_1=\bar{W}_2'\bar{W}_2=I_{\min(m,T)\times\min(m,T)}$.
Also
\[
\hat{A}^{(2)}=U\pmatrix{
0 \!& 0\cr
0 \!& W_1\Lambda(\tilde{B}_{22})W_2'
}V'=U\bar{W}_1
\pmatrix{
0 \!& 0\cr
0 \!&
\Lambda(\tilde{B}_{22})
}\bar{W}_2'V'=:U\bar{W}_1D\bar{W}_2'V'.
\]
We now represent $\hat{A}^{(2)}$ as a finite sum of matrices
$\hat{A}^{(2)}=\sum_{j=1}^{R'}\hat{A}^{(2)}_j$ with
\[
\hat{A}^{(2)}_i=U\bar{W}_1D_i\bar{W}_2'V'
\]
and
\[
D_i=\pmatrix{ 0 & 0\cr0 & \Lambda_i
},
\]
where the $r'\times r'$ diagonal matrix $\Lambda_i$ has the form
$\Lambda_i=\operatorname{diag}(\lambda_j I_{\{j\in I_i\}})$, $i\geq1$. We
denote here by $I_1$ the set of $ar$ indices from
$\{1,\dots,\min(m,T)\}$ corresponding to the $ar$ largest in absolute
value diagonal entries of $\Lambda$, by $I_2$ the set of indices
corresponding to the next $ar$ largest in absolute value diagonal
entries $\lambda_j$, etc. Clearly, the matrices $\hat{A}^{(2)}_k$ are
mutually orthogonal:
$\operatorname{tr}((\hat{A}^{(2)}_j)'\hat{A}^{(2)}_k)=0$ for $j\not
= k$
and $\operatorname{rank}(\hat{A}^{(2)}_j)\leq ar$. Moreover,
$\hat{A}^{(2)}_i$ is orthogonal to $\hat{A}^{(1)}-A^*$.

Let $\sigma_1\geq\sigma_2\geq\cdots$ be the singular values of
$\hat{A}^{(2)}$, then $\sigma_1\geq\cdots\geq\sigma_{ar}$ are the
singular values of $\hat{A}^{(2)}_1$, $\sigma_{ar+1}\geq\cdots\geq
\sigma_{2ar}$ those of $\hat{A}^{(2)}_2$, etc. By construction, we have
$\operatorname{Card}(I_i)=ar$ for all $i$, and for all $k\in I_{i+1}$
\[
\sigma_k\leq\min_{j\in I_i}\sigma_j\leq\biggl(\frac{1}{ar}\sum_{j\in
I_i}\sigma_j^p\biggr)^{1/p}.
\]
Thus,
\[
\sum_{k\in I_{i+1}}\sigma_k^2\leq ar\biggl(\frac{1}{ar}\sum_{j\in
I_i}\sigma_j^p\biggr)^{2/p}
\]
from which one can deduce for all $j\geq2$:
\[
\bigl\Vert\hat{A}_j^{(2)}\bigr\Vert_{S_2}=\biggl(\sum_{k\in
I_{j}}\sigma_k^2\biggr)^{1/2}\leq(ar)^{1/2-1/p}\biggl(\sum
_{k\in
I_{j-1}}\sigma_k^p\biggr)^{1/p}=(ar)^{1/2-1/p}\bigl\Vert
\hat{A}_{j-1}^{(2)}\bigr\Vert_{S_p}
\]
and consequently
\[
\sum_{j\geq2}\bigl\Vert
\hat{A}_j^{(2)}\bigr\Vert_{S_2}\leq(ar)^{1/2-1/p}\sum_{j\geq
1}\bigl\Vert\hat{A}_j^{(2)}\bigr\Vert_{S_p}.
\]
Because of the elementary inequality $x^{1/p}+y^{1/p}\leq(x+y)^{1/p}$
for any nonnegative $x,y$ and $0<p\leq1$,
\begin{eqnarray*}
\sum_{j\geq2}\bigl\Vert\hat{A}_j^{(2)}\bigr\Vert_{S_p} &=&\sum
_{j\geq
2}\biggl(\sum_{k\in I_j}\sigma_k^p\biggr)^{1/p}\leq\biggl(\sum_{j\geq2}\sum_{k\in
I_j}\sigma_k^p\biggr)^{1/p}\\&\leq&\biggl(\sum_k\sigma_k^p\biggr)^{1/p}=\bigl\Vert\hat
{A}^{(2)}\bigr\Vert_{S_p}.
\end{eqnarray*}
Therefore,
\begin{eqnarray*}
\sum_{j\geq2}\bigl\Vert\hat{A}_j^{(2)}\bigr\Vert_{S_2}&\leq&
(ar)^{1/2-1/p}\bigl\Vert\hat{A}^{(2)}\bigr\Vert
_{S_p}\\
&\leq&3^{1/p}(ar)^{1/2-1/p}\bigl\Vert\hat
{A}^{(1)}-A^*\bigr\Vert_{S_p}\qquad\mbox{[using inequality (\ref{eq:
1--2})]}\\
&\leq&3^{1/p}(ar)^{1/2-1/p}
(2r)^{1/p-1/2}\bigl\Vert
\hat{A}^{(1)}-A^*\bigr\Vert_{S_2},
\end{eqnarray*}
where the last inequality results from
$\operatorname{rank}(\hat{A}^{(1)}-A^*)\leq2r$ and
\[
\biggl(\frac{1}{2r}\sum_{k\leq
2r}\sigma_k^p\biggr)^{1/p}\leq\biggl(\frac{1}{2r}\sum_{k\leq2r}\sigma_k^2\biggr)^{1/2}.
\]
Finally,
%
\begin{equation}\label{sc1}
\sum_{j\geq
2}\bigl\Vert\hat{A}_j^{(2)}\bigr\Vert_{S_2}\leq3^{1/p}\biggl(\frac
{a}{2}\biggr)^{1/2-1/p}
\bigl\Vert\hat{A}^{(1)}-A^*\bigr\Vert_{S_2}.
\end{equation}
We now proceed with the final argument. First, note that
$\operatorname{rank}((\hat{A}^{(1)}-A^*)+\hat{A}^{(2)}_1)\leq(2+a)r$. Next,
by the triangular inequality, the restricted isometry condition and the
orthogonality of $\hat{A}^{(2)}_j$ and $\hat{A}^{(1)}-A^*$ we obtain
\begin{eqnarray}\label{sc}
\nu\hat{d}_{2,N}(\hat{A},A^*)
&=&
\nu\vert\LL(\hat{A}-A^*)\vert_2\nonumber
\\
&\geq&
\nu\bigl\vert\LL\bigl(\hat{A}^{(1)}-A^*+\hat{A}^{(2)}_1\bigr)\bigr\vert_2- \nu\sum_{j\geq2}\bigl\vert\LL\bigl(\hat{A}^{(2)}_j\bigr)\bigr\vert_2
\\
&\geq&
\bigl(1-\delta_{(2+a)r}\bigr)\bigl\Vert\hat{A}^{(1)}-A^*+\hat{A}^{(2)}_1\bigr\Vert_{S_2}-(1+\delta_{ar})\sum_{j\geq2}\bigl\Vert\hat{A}^{(2)}_j\bigr\Vert_{S_2}\nonumber
\\
&\geq&
\bigl\Vert\hat{A}^{(1)}-A^*\bigr\Vert_{S_2}\biggl(\bigl(1-\delta_{(2+a)r}\bigr)-(1+\delta_{ar})3^{1/p}\biggl(\frac{a}{2}\biggr)^{1/2-1/p}\biggr).\nonumber
\end{eqnarray}
Define
\[
a=a(p)=\min\bigl\{k\in\N\dvtx k>\bigl(6^{1/p}/\sqrt{2}\bigr)^{2p/(2-p)}\bigr\}.
\]
Then $1-3^{1/p}(a/2)^{1/2-1/p}>0$. Now,
$\delta_{(2+a)r}\geq\delta_{ar}$, and thus
\[
\bigl(1-\delta_{(2+a)r}\bigr)-(1+\delta_{ar})3^{1/p}\biggl(\!\frac{a}{2}\!\biggr)^{\!1/2-1/p}
\geq\biggl(\!1-3^{1/p}\biggl(\frac{a}{2}\biggr)^{\!1/2-1/p}\biggr)-2\delta_{(2+a)r}>0
\]
whenever
%
\begin{equation}\label{eq: hallo-1}
\delta_{(2+a)r}<\frac{1}{2}\biggl(1-3^{1/p}\biggl(\frac{a}{2}\biggr)^{1/2-1/p}\biggr).
\end{equation}
In case of (\ref{eq: hallo-1}), there exists a universal constant
$\kappa=\kappa(p)$ such that
%
\begin{equation}\label{eq: RIb}
\nu^2\hat{d}_{2,N}(\hat{A},A^*)^2\geq\kappa\bigl\Vert
\hat{A}^{(1)}-A^*\bigr\Vert_{S_2}^2.
\end{equation}
Now, inequalities (\ref{eq: ub}) and (\ref{eq: RIb}) yield
%
\begin{equation}\label{s0}
\quad\quad\kappa\bigl\Vert\hat{A}^{(1)}-A^*\bigr\Vert_{S_2}^2\leq6\tau\nu^2
\bigl\Vert\hat{A}^{(1)}-A^*\bigr\Vert_{S_p}^p\leq6\tau\nu^2
(2r)^{1-p/2}\bigl\Vert\hat{A}^{(1)}-A^*\bigr\Vert_{S_2}^p,\hspace*{-15pt}
\end{equation}
where the second inequality results from the fact that
$\operatorname{rank}(\hat{A}^{(1)}-A^*)\leq2r$, which implies
%
\begin{equation}\label{s1}
\bigl\Vert
\hat{A}^{(1)}-A^*\bigr\Vert_{S_p}\leq(2r)^{1/p-1/2}\bigl\Vert
\hat{A}^{(1)}-A^*\bigr\Vert_{S_2}.
\end{equation}
From (\ref{s0}), we obtain
%
\begin{equation}\label{s2}
\kappa\bigl\Vert
\hat{A}^{(1)}-A^*\bigr\Vert_{S_2}^{2-p}\leq6\tau\nu^2 (2r)^{1-p/2}.
\end{equation}
Furthermore, from (\ref{eq: ub}), (\ref{s1}) and (\ref{s2}) we find
\begin{eqnarray}\label{s3}
\hat{d}_{2,N}(\hat{A},A^*)^2&\le&6\tau(2r)^{1-p/2} \bigl\Vert
\hat{A}^{(1)}-A^*\bigr\Vert_{S_2}^p \nonumber\\[-8pt]\\[-8pt]
&\le&2 r (6\tau)^{2/(2-p)}
\kappa^{-p/(2-p)}\nu^{2p/(2-p)}.\nonumber
\end{eqnarray}
This proves (\ref{eq:mainn}). It remains to prove (\ref{eq:mainn1}). We
first demonstrate (\ref{eq:mainn1}) for $q=2$, then for $q=p$, and
finally obtain (\ref{eq:mainn1}) for all $q\in[p,2]$ by Schatten norm
interpolation.

Using (\ref{sc1}), (\ref{sc}), (\ref{s3}), we find
\begin{eqnarray*}
\bigl(1-\delta_{(2+a)r}\bigr)\bigl\Vert
\hat{A}^{(1)}-A^*+\hat{A}^{(2)}_1\bigr\Vert_{S_2} & \le&
\nu\hat{d}_{2,N}(\hat{A},A^*) + (1+\delta_{ar})\sum_{j\geq
2}\bigl\Vert\hat{A}^{(2)}_j\bigr\Vert_{S_2}
\\
&\le& C\sqrt{r}\tau^{1/(2-p)} \nu^{2/(2-p)}
\end{eqnarray*}
for some constant $C=C(p)>0$. This and again (\ref{sc1}) yield
\[
\Vert\hat{A}-A^*\Vert_{S_2} \le\bigl\Vert
\hat{A}^{(1)}-A^*+\hat{A}^{(2)}_1\bigr\Vert_{S_2} + \sum_{j\geq
2}\bigl\Vert\hat{A}^{(2)}_j\bigr\Vert_{S_2} \le
C'\sqrt{r}\tau^{1/(2-p)} \nu^{2/(2-p)}
\]
for some constant $C'=C'(p)>0$. Thus, we have proved (\ref{eq:mainn1})
for $q=2$. Next, using inequalities (\ref{triang}) and (\ref{eq: 1--2})
we obtain
\[
\Vert\hat{A}-A^*\Vert_{S_p}^p\le
\bigl\Vert\hat{A}^{(1)}-A^*\bigr\Vert_{S_p}^p +
\bigl\Vert\hat{A}^{(2)}\bigl\Vert_{S_p}^p \le
4\bigl\Vert\hat{A}^{(1)}-A^*\bigr\Vert_{S_p}^p.
\]
Combining this with (\ref{s1}) and (\ref{s2}) we get (\ref{eq:mainn1})
for $q=p$. Finally, (\ref{eq:mainn1}) for arbitrary $q\in[p,2]$ follows
from the norm interpolation formula
\[
\Vert A \Vert_{S_q}^q \le\Vert A
\Vert_{S_p}^{p(2-q)/(2-p)} \Vert A
\Vert_{S_2}^{2(q-p)/(2-p)};
\]
cf. Lemma \ref{lemma7a} of Section \ref{sec: embeddings} with
$\theta=\frac{p(2-p)}{q(2-q)}$.
\end{pf*}

\section{Proofs of the lemmas}
\label{sec:proof_lem}
\mbox{}
\begin{pf*}{Proof of Lemma \ref{lemma:bernstein}}
First, observe that
\[
\Vert\mathbf{M} \Vert_{S_{\infty}} = \mathop{\sup_{u\in
\R^T:}}_{\vert u\vert_2=1}\vert\mathbf{M}u\vert_2\leq\sqrt{m}\max_{1\leq j\leq
m}\mathop{\sup_{u\in\R^T:}}_{\vert u\vert
_2=1}\vert
u'\bar{\eta}_j\vert,
\]
with vectors
$\bar{\eta}_j=N^{-1}\sum_{i=1}^N\xi_i X_{i(j,\cdot)}$. Consequently,
for any $t>0$,
\begin{eqnarray*}
\PP\Biggl(\Vert\mathbf{M}\Vert_{S_{\infty}}\geq t\sqrt{\frac {m\log
m}{N}}\Biggr)&\leq&\PP\Biggl(\sqrt{m}\max_{1\leq j\leq m}\vert
\bar{\eta}_j\vert_2\geq t\sqrt{\frac{m\log m}{N}}\Biggr)\\&\leq&
m\max_{1\leq j\leq m}\PP\Biggl(\vert\bar{\eta}_j\vert_2\geq
t\sqrt{\frac{\log m}{N}}\Biggr).
\end{eqnarray*}
To proceed with the evaluation of the latter probability, we use the
following concentration bound [Pinelis and Sakhanenko (\citeyear{ps85})].

\begin{lemma}\label{lemma: Hilbert}
Let $\zeta_1,\ldots,\zeta_N$ be independent zero mean random variables
in a separable Hilbert space $\HH$ such that
%
\begin{equation}\label{bh}
\sum_{i=1}^N\E\Vert\zeta_i\Vert_{\HH}^l\leq\frac
{1}{2}l!B^2L^{l-2},\qquad l=
2,3,\ldots,
\end{equation}
with some finite constants $B,L> 0$. Then
\[
\PP\Biggl(\Biggl\Vert\sum_{i=1}^N\zeta_i\Biggr\Vert_{\HH}\geq
x\Biggr)\leq2\exp\biggl(-\frac{x^2}{2B^2+2xL}\biggr)\qquad \forall x>0.
\]
\end{lemma}

Setting $\zeta_i=\xi_iX_{i(j,\cdot)}$, $\HH=\R^T$, note first
that, by the Bernstein condition~(\ref{eq: Bernstein}),
\begin{eqnarray*}
\sum_{i=1}^N\E\Vert\zeta_i\Vert_{\HH}^l&=&\E\vert
\xi_i\vert^l\sum_{i=1}^N\bigl\vert X_{i(j,\cdot)}\bigr\vert
_2^l\\
&\leq&\frac{1}{2}l!\sigma^2H^{l-2}\Biggl(\max_j\sum_{i=1}^N\bigl\vert
X_{i(j,\cdot)}\bigr\vert_2^2\Biggr)\max_{i,j}\bigl\vert X_{i(j,\cdot
)}\bigr\vert_2^{l-2}\\
&\leq&\frac{1}{2}l!B^2L^{l-2},
\end{eqnarray*}
where $B^2=\sigma^2S_{\mathrm{row}}^2N$ and $L= H_{\mathrm{row}}H$, that is, condition
(\ref{bh}) is satisfied. Now an application of Lemma \ref{lemma:
Hilbert} yields for any $t>0$
\begin{eqnarray*}
\PP\Biggl(\vert\bar{\eta}_j\vert_2\geq t\sqrt{\frac{\log
m}{N}}\Biggr)&=&\PP\Biggl(\Biggl\vert\frac{1}{\sqrt{N}}\sum_{i=1}^N\xi
_i\Biggr\vert_2>t\sqrt{\log m}\Biggr)\\
&\leq&2\exp\biggl(-\frac{N(\log m)t^2}{2B^2+2tL\sqrt{N\log m}}\biggr)\\
&=&2\exp\biggl(-\frac{N(\log m)t^2}{2\sigma^2S_{\mathrm{row}}^2N+2tL\sqrt{N\log m}}\biggr).
\end{eqnarray*}
Define $ t = \sqrt{2D\sigma^2S_{\mathrm{row}}^2}+2DL\sqrt{\frac{\log m}{N}} $
for some $D>1$. Then
\[
\frac{t^2}{\bar{B}+\bar{L}t}\geq
D,\qquad\mbox{where }\bar{B}=2\sigma^2S_{\mathrm{row}}^2,\bar{L}=2L\sqrt{\frac
{\log
m}{N}}.
\]
With this choice of $t$,
\[
\PP\Biggl(\vert\bar{\eta}_j\vert_2\geq t\sqrt{\frac{\log
m}{N}}\Biggr)\leq 2\exp(-D\log m)=2m^{-D}
\]
and therefore $ \PP(\Vert\mathbf{M}\Vert_{S_{\infty}}\geq\tau_{\mathrm{row}}) \leq2 m^{1-D} $, where
\[
\tau_{\mathrm{row}}=\Biggl(\sqrt{2D\sigma^2S_{\mathrm{row}}^2}+2DH_{\mathrm{row}}H \sqrt{\frac
{\log
m}{N}}\Biggr)\sqrt{\frac{m\log m}{N}}.
\]
Similarly, using $\Vert\mathbf{M}\Vert_{S_{\infty}}=\sup_{\vert v\vert _2=1}\vert
v'\mathbf{M}\vert_2$, and assuming (\ref{eq: con 1a}) and (\ref{eq:
con 1b}), we get $ \PP(\Vert\mathbf{M}\Vert_{S_{\infty}}\geq\tau_{\mathrm{col}}) \leq2T^{1-D}$, where
\[
\tau_{\mathrm{col}}=\Biggl(\sqrt{2D\sigma^2S_{\mathrm{col}}^2}+2DH_{\mathrm{col}}H\sqrt{\frac{\log
T}{N}}\Biggr)\sqrt{\frac{T\log T}{N}}.
\]
\upqed\end{pf*}

\begin{pf*}{Proof of Lemma \ref{lemma:stochterm_matr_compl}}
The matrix $\mathbf{M}=\frac{1}{N}\sum_{i=1}^N\xi_iX_i$ is a sum of i.i.d.
random matrices. Therefore, part (ii) of the lemma follows by direct
application of the large deviations inequality of Nemirovski (\citeyear{nem04}).

To prove part (i) of the lemma, we use bounds on maximal eigenvalues of
subgaussian matrices due to Mendelson, Pajor and Tomczak-Jaegermann (\citeyear{mptffjf07}); see also
Vershynin (\citeyear{vershf07}). However, direct application of these bounds (based on
the overall subgaussianity) does not lead to rates that are accurate
enough for our purposes. We therefore need to refine the argument using
the specific structure of the matrices. Note first that
\[
\| \mathbf{M} \|_{S_{\infty}}= \max_{v\in\mathcal{S}^{T-1}} |\mathbf{M}v|_2 =
\max_{u\in\mathcal{S}^{m-1},v\in\mathcal{S}^{T-1}} u'\mathbf{M}v,
\]
where $\mathcal{S}^{m-1}$ is the unit sphere in $\R^m$. Therefore,
denoting by $\mathcal{M}_m$ and $\mathcal{M}_T$ the minimal $1/2$-nets in
Euclidean metric on $\mathcal{S}^{m-1}$ and $\mathcal{S}^{T-1}$, respectively,
we easily get
\[
\| \mathbf{M} \|_{S_{\infty}}\le2\max_{v\in\mathcal{M}_T} |\mathbf{M}v|_2
\le
4 \max_{u\in\mathcal{M}_m,v\in\mathcal{M}_T} |u'\mathbf{M}v|.
\]
Now, $\operatorname{Card}(\mathcal{M}_m)\le5^m$ [cf. Kolmogorov and
Tikhomirov (\citeyear{ktff59})] so that by the union bound, for any $\tau>0$,
%
\begin{equation}\label{mcc1}
\PP(\Vert\mathbf{M}\Vert_{S_{\infty}}\geq\tau)
\leq5^{m+T}\max_{u\in\mathcal{M}_m,v\in\mathcal{M}_T} \PP(|u'\mathbf{M}v| \ge
\tau/4).
\end{equation}
It remains to bound the last probability in (\ref{mcc1}) for fixed
$u,v$. Let us fix some $u\in\mathcal{S}^{m-1},v\in\mathcal{S}^{T-1}$ and
introduce the random event
\[
\mathcal{A}= \Biggl\{\frac1{N}\sum_{i=1}^N (u'X_i v)^2 \le\frac{5(m+T)}{N}\Biggr\}.
\]
Note that $\E(u'X_i v)^2= \sum_{k=1}^m\sum_{l=1}^T u_k^2v_l^2
\PP(X_1=e_k(m)e_l'(T))= (mT)^{-1}|u|_2^2 \times|v|_2^2=(mT)^{-1}$, and
consider the zero-mean random variables $\eta_i = (u'X_i v)^2 -
\E(u'X_i v)^2=(u'X_i v)^2- (mT)^{-1}$. We have \mbox{$|\eta_i|\!\le\!2 \max_i
(u'X_i v)^2 \!\le\!2 |u|_2^2 |v|_2^2=2$}. Furthermore,
\begin{eqnarray*}
\E(\eta_i^2) &\le&\E(u'X_i v)^4 \le\sum_{k=1}^m\sum_{l=1}^T u_k^4v_l^4
\PP\bigl(X_1=e_k(m)e_l'(T)\bigr)\\
& =& (mT)^{-1}\sum_{k=1}^m u_k^4 \sum_{l=1}^T
v_l^4 \le(mT)^{-1}.
\end{eqnarray*}
Therefore, using Bernstein's inequality and the condition $(m+T)/N
> (mT)^{-1}$ we get
\begin{eqnarray}\label{mcc2}
\PP(\mathcal{A}^c) &\le&2 \exp\biggl(-\frac{N(4(m+T)/N)^2}{2(mT)^{-1}
+(4/3)(4(m+T)/N)}\biggr)\nonumber\\[-8pt]\\[-8pt]
&\le&2 \exp\bigl(-2(m+T)\bigr),\nonumber
\end{eqnarray}
where $\mathcal{A}^c$ is the complement of $\mathcal{A}$. We now bound the
conditional probability
\[
\PP(|u'\mathbf{M}v| \ge\tau/4|X_1,\dots,X_N) = \PP\Biggl(\Biggl|\frac1{N}\sum_{i=1}^N
\xi_i(u'X_iv)\Biggr| \ge\tau/4\bigl|X_1,\dots,X_N\Biggr).
\]
Note that conditionally on $X_1,\dots,X_N$, the $\xi_i(u'X_iv)$ are
independent zero-mean random variables with
\[
\sum_{i=1}^N\E(|\xi_i(u'X_iv)|^l|X_1,\ldots,X_N)
\le\E|\xi_1|^l \sum_{i=1}^N|u'X_iv|^2  \qquad\forall l\ge2,
\]
where we used the fact that $|u'X_iv|^{l-2} \le(|u|_2 |v|_2)^{l-2}= 1$
for $l\ge2$. This and the Bernstein condition (\ref{eq: Bernstein})
yield that, for $(X_1,\dots,X_N)\in\mathcal{A}$,
\[
\sum_{i=1}^N\E(|\xi_i(u'X_iv)|^l|X_1,\dots,X_N) \le\frac{l!}{2}B^2
H^{l-2}
\]
with $B^2=5(m+T)\sigma^2$. Therefore, by Lemma \ref{lemma: Hilbert},
for $(X_1,\dots,X_N)\in\mathcal{A}$ we have
%
\begin{equation}\label{mcc3}
\quad\quad\PP(|u'\mathbf{M}v| \ge\tau/4|X_1,\dots,X_N) \le
2\exp\biggl(-\frac{N^2\tau^2/16}{10\sigma^2(m+T)+N\tau H/2}\biggr).\hspace*{-10pt}
\end{equation}
For $\tau$ defined in (\ref{eq:stochterm0}) the last expression does
not exceed $2\exp(-D(m+T))$. Together with (\ref{mcc1}) and
(\ref{mcc2}), this proves the lemma.
\end{pf*}


\begin{pf*}{Proof of Lemma \ref{lemma:stochterm_matr_compl_1}}
We act as in the proof of Lemma \ref{lemma:stochterm_matr_compl} but
since the matrices $X_i$ are now deterministic, we do not need to
introduce the event $\mathcal{A}$. By the definition of $\phi_{\max}(1)$,
\[
\frac1{N}\sum_{i=1}^N (u'X_i v)^2 = |\mathcal{L}(uv')|_2^2 \le
\phi_{\max}^2(1)\|uv'\|_{S_2}^2 = \phi_{\max}^2(1)
\]
for all $u\in\mathcal{S}^{m-1},v\in\mathcal{S}^{T-1}$. Hence,
$\frac1{N}\sum_{i=1}^N \xi_i(u'X_iv)$ is a zero-mean Gaussian random
variable with variance not larger than $\phi_{\max}^2(1)\sigma^2/N$.
Therefore,
\[
\PP(|u'\mathbf{M}v| \ge\tau/4) \le 2\exp\biggl(-\frac{N\tau^2}{32
\phi_{\max}^2(1)\sigma^2}\biggr).
\]
For $\tau$ as in (\ref{eq:stochterm00}) the last expression does not
exceed $2\exp(-D(m+T))$. Combining this with (\ref{mcc1}), we get the
lemma.
\end{pf*}


\begin{pf*}{Proof of Lemma \ref{lemma:stochterm_matr_compl_2}}
We proceed again as in the proof of Lemmas~\ref{lemma:stochterm_matr_compl}
and~\ref{lemma:stochterm_matr_compl_1}. Denote by $\Omega$ the set of
pairs $(k,l)$ such that $\{X_1,\dots,X_N\}=\{e_k(m)e_l'(T),\break (k,l)\in
\Omega\}$ (recall that all $X_i$ are distinct by assumption). Then
%
\begin{equation}\label{mcc5}
\sum_{i=1}^N (u'X_iv)^2=\sum_{(k,l)\in\Omega}u_k^2v_l^2 \le|u|_2^2
|v|_2^2 =1
\end{equation}
for any $u\in\mathcal{S}^{m-1},v\in\mathcal{S}^{T-1}$. Hence, under the
assumptions of part (i) of the lemma,
\[
\PP(|u'\mathbf{M}v| \ge\tau/4) \le 2\exp\biggl(-\frac{N^2\tau ^2}{32\sigma^2}\biggr)
\]
which does not exceed $2\exp(-D(m+T))$ for $\tau$ defined in
(\ref{eq:stochterm010}). Combining this with (\ref{mcc1}) we get part
(i) of the lemma. To prove part (ii) we note that, as in the proof of
Lemma \ref{lemma:stochterm_matr_compl}, $|u'X_iv|^{l-2} \le1$ for
$l\ge2$. This and (\ref{mcc5}) yield
\[
\sum_{i=1}^N\E(|\xi_i(u'X_iv)|^l) \le\frac{l!}{2}B^2 H^{l-2}\qquad
\forall
l\ge2,
\]
with $B^2=\sigma^2$. Therefore, by Lemma \ref{lemma: Hilbert}, we have
\[
\PP(|u'\mathbf{M}v| \ge\tau/4) \le
2\exp\biggl(-\frac{N^2\tau^2/16}{2\sigma^2+N\tau H/2}\biggr),
\]
and we complete the proof of (ii) in the same way as in Lemmas
\ref{lemma:stochterm_matr_compl} and
\ref{lemma:stochterm_matr_compl_1}.

Part (iii) follows by an application of Theorem 2.1, Tropp (\citeyear{troppgg10}),
after replacing every $X_i$ by its self-adjoint dilation [see Paulsen
(\citeyear{paulsenf86})].
\end{pf*}


For the proof of Lemma \ref{nonconv} we will need some notation. The
$p$th Schatten class of $M\times M$-matrices is denoted by $S_p^M$, and
we write
\[
\BB(S_p^M)=\{A\in\R^{M\times M}\dvtx \Vert A\Vert_{S_p}\leq
1\}
\]
for the corresponding closed Schatten-$p$ unit ball in $\R^{M\times
M}$. For any pseudo-metric space $(\TT, d)$ and any $\varepsilon>0$, we
define the covering number
\[
\NN(\TT, d,\varepsilon)=\min\Bigl\{\operatorname{Card}(\TT_0)\dvtx \TT_0\subset\TT\mbox{ and }\inf_{s\in\TT_0}d(t,s)\leq\varepsilon
\mbox{ for
all }t\in\TT\Bigr\}.
\]
In other words, $\NN(\TT, d,\varepsilon)$ is the smallest number of
closed balls of radius $\varepsilon$ in the metric $d$ needed to cover
the set $\TT$. We will sometimes write $\NN(\TT,
\|\cdot\|,\varepsilon)$ instead of $\NN(\TT, d,\varepsilon)$ if the
metric $d$ is associated with the norm $\|\cdot\|$. The empirical norm
$\Vert\cdot\Vert_{2,N}$ corresponds to $\hat{d}_{2,N}$, that
is, for all $A\in\R^{M\times M}$,
\[
\Vert A\Vert_{2,N}^2=\frac{1}{N} \sum_{j=1}^N\trace(A'X_j)^2.
\]

\begin{pf*}{Proof of Lemma \ref{nonconv}}
Let us first assume that $m=T\equiv M$. Since
\[
\sup_{B\in\R^{M\times M}}\biggl\vert
\frac{(1/\sqrt{N})\sum_{i=1}^N\xi_i\operatorname
{tr}(B'X_i)}{\Vert
B\Vert_{2,N}^{1-p/(2-p)}\Vert
B\Vert_{S_p}^{p/(2-p)}}\biggr\vert=\sup_{B\in\BB
(S_p^M)}\biggl\vert\frac{(1/\sqrt{N})\sum_{i=1}^N\xi
_i\operatorname{tr}(B'X_i)}{\Vert
B\Vert_{2,N}^{1-p/(2-p)}}\biggr\vert,
\]
the expression on the LHS of (\ref{eq:nonconv}) is not greater than
\begin{eqnarray*}
&&\frac{\sqrt{M}}{\sqrt{p}\sqrt{N}}\hat{d}_{2,N}(\hat
{A},A^*)^{1-p/(2-p)}\Vert\hat{A}-A^*\Vert
_{S_p}^{p/(2-p)}\\
&&\qquad{}\times
\sup_{B\in\BB(S_p^M)}\biggl\vert\frac
{(M/p)^{(p-2)/(2p)}N^{-1/2}\sum_{i=1}^N\xi_i\operatorname
{tr}(B'X_i)}{((M/p)^{(p-2)/(2p)}\Vert
B\Vert_{2,N})^{1-p/(2-p)}}\biggr\vert.
\end{eqnarray*}
Due to the linear dependence in $M$ of the $\varepsilon$-entropies of
the quasi-convex Schatten class embeddings $ S_p^M\hookrightarrow S_2^M
$ (cf. Corollary \ref{cor: entropy}) and the fact that the required
bound should be uniform in $M$ and in $p$ for $p\searrow0$, we
introduced an additional weighting by $(M/p)^{(p-2)/2p}$. Now
define
\[
\GG_{M,p}=\bigl\{A\in\R^{M\times M}\dvtx
(M/p)^{(2-p)/(2p)}A\in\BB(S_p^M)\bigr\}.
\]
By the entropy bound of Corollary \ref{cor: entropy} and the uniform
boundedness condition~(\ref{eq: UB}),
\[
\log\NN(\GG_{M,p},
\hat{d}_{2,N},\varepsilon)\leq\log\NN\bigl(\GG_{M,p},\sqrt{c_0}\|
\cdot\|_{S_2},\varepsilon\bigr)\leq
p\alpha_0(p)\bigl(\varepsilon/\sqrt{c_0}\bigr)^{-2p/(2-p)}
\]
whence
%
\begin{equation}\label{eq: entropy integral}
\quad\quad\int_0^{\delta}\sqrt{\log\NN(\GG_{M,p}, \hat
{d}_{2,N},\varepsilon)}
\,d\varepsilon\leq
c_0^{p/(2(2-p))}p\alpha_0(p)\frac{2-p}{2-2p}\delta^{1-p/(2-p)}.\hspace*{-10pt}
\end{equation}
We remark that due to the order specification of $\alpha_0$ in
Corollary \ref{cor: entropy}, the expression
%
\begin{equation}
c_0^{p/(2(2-p))}p\alpha_0(p)\frac{2-p}{2-2p}\label{eq:
uniformity}
\end{equation}
is uniformly bounded as long as $p$ stays uniformly bounded away from
$1$. Note that for $p=1$ the entropy integral on the LHS in (\ref{eq:
entropy integral}) does not converge.

\begin{cla}
For any $q\in(0,1)$, there exist
constants $c(q)$ and $c'(q)$, such that for all $0<p\leq q$, all
$0<\delta\leq\sqrt{c_0}$ and uniformly in $M$ and $N$,
%
\begin{equation}
\PP\Biggl(\mathop{\sup_{B\in\GG_{M,p}:}}_{\Vert B\Vert _{2,N}\leq
\delta}\Biggl\vert\frac{1}{\sqrt{N}}\sum_{j=1}^N\xi
_j\operatorname{tr}(X_j'B)\Biggr\vert\geq T\Biggr)\leq
c(q)\exp\biggl(-\frac{T^2}{c(q)^2\delta^2}\biggr)
\end{equation}
for
all $T\geq c'(q)\delta^{1-p/(2-p)}$.
\end{cla}

\begin{pf}
The bound is essentially stated
in van de
Geer (\citeyear{vdgeerf00}) as Lemma 3.2 [further referred to as VG(00)]. The constant
in VG(00) depends neither on the $\Vert
\cdot\Vert_{2,N}$-diameter of the function class nor on the
function class itself and is valid, in particular, for $\varepsilon=0$,
in the notation of VG(00). The uniformity in $0<p\leq q$ follows from
the uniform boundedness of (\ref{eq: uniformity}) over $p\in(0,q]$. The
required case corresponds to $K=\infty$ in the notation of VG(00). Its
proof follows by taking $\varepsilon=0$ and applying the theorem of
monotone convergence as $K\rightarrow\infty$, since the RHS of the
inequality is independent of $K$.
\end{pf}

\begin{cla}
For any $q\in(0,1)$, there exists a
constant $C(q)$ such that for any $0<p\leq q$
%
\begin{equation}\label{eq: eb}
\quad\quad\PP\biggl(\sup_{B\in\GG_{M,p}}\biggl\vert\frac{(1/\sqrt{N})\sum
_{j=1}^N\xi_j\operatorname{tr}(B'X_j)}{\Vert
B\Vert_{2,N}^{1-p/(2-p)}}\biggr\vert\geq T\biggr) \leq
C(q)\exp\bigl(-T^2M/C(q)^2\bigr)\hspace*{-15pt}
\end{equation}
for all $T\geq C(q)$.
\end{cla}

\begin{pf}
First, observe that
\begin{eqnarray*}
\sup_{A\in\GG_{M,p}}\Vert A\Vert_{2,N} &\leq&
\sqrt{c_0}\sup_{A\in\GG_{M,p}}\Vert A\Vert_{S_2}
\\
&\leq&
\sqrt{c_0}(M/p)^{(p-2)/(2p)}\sup_{A\in\BB(S_2^M)}\Vert
A\Vert_{S_2}=\sqrt{c}_0(M/p)^{(p-2)/(2p)},
\end{eqnarray*}
where the last inequality follows from $\BB(S_p^M)\subset\BB(S_2^M)$.
Define the decomposition of $\GG_{M,p}$
%
\begin{eqnarray}
&&\GG_{M,p}^{(k)}=\bigl\{A\in\GG_{M,p}\dvtx
(1/2)^{k}\sqrt{c_0}(M/p)^{(p-2)/(2p)}\nonumber
\\
&&\hphantom{\GG_{M,p}^{(k)}=\bigl\{A\in\GG_{M,p}\dvtx}\qquad\leq
\Vert
A\Vert_{2,N}\leq(1/2)^{k-1}\sqrt{c_0}(M/p)^{(p-2)/(2p)}\bigr\},\nonumber
\\
&&\eqntext{k\in\N.}
\end{eqnarray}
Then by peeling-off the class $\GG_{M,p}$, we obtain together with
claim I
for all $T\geq c'(q)$
\begin{eqnarray}\label{eq: c-0}
&&\PP\biggl(\sup_{B\in\GG_{M,p}}\biggl\vert\frac{(1/\sqrt{N})\sum
_{j=1}^N\xi_j\operatorname{tr}(B'X_j)}{\Vert B\Vert
_{2,N}^{1-p/(2-p)}}\biggr\vert\geq T\biggr)\nonumber
\\
&&\qquad\leq
\sum_{k=1}^{\infty}\PP\Biggl( \sup_{B\in\GG_{M,p}^{(k)}}\Biggl\vert\frac{1}{\sqrt{N}}\sum_{j=1}^N\xi_j\operatorname{tr}(B'X_j)\Biggr\vert\nonumber
\\
&&\qquad\hphantom{\qquad\leq
\sum_{k=1}^{\infty}\PP\Biggl(}\geq T\bigl((1/2)^{k}\sqrt{c_0}(M/p)^{(p-2)/(2p)}\bigr)^{1-p/(2-p)}\Biggr)\nonumber
\\
&&\qquad\leq
\sum_{k=1}^{\infty}c(q)\exp\biggl(-\frac{T^{2}(1/2)^{2}((1/2)^{k}\sqrt{c_0}(M/p)^{(p-2)/(2p)})^{-2p/(2-p)}}{c(q)^2}\biggr)\nonumber
\\
&&\qquad\leq\sum_{k=1}^{\infty}c(q)\exp\biggl(-\frac{T^2M2^{k(2p)/(2-p)}C_0(q)}{4pc(q)^2}\biggr)
\end{eqnarray}
with the definition
\[
C_0(q)= \inf_{0<p\leq q}c_0^{-p/(2-p)}.
\]
It remains to note that the last sum in (\ref{eq: c-0}) is bounded by $
C(q)\exp(-T^2M/\break C(q)^2) $ uniformly in $0<p\leq q$ whenever $T\geq C(q)$ for some suitable
constant $C(q)$. This follows from the fact that
\[
\sum_{k=1}^{\infty}\exp\bigl(-p^{-1}2^{k(2p)/(2-p)}\bigr)\leq
\sum_{k=1}^{\infty}\frac{1}{p^{-1}2^{k(2p)/(2-p)}+1}\leq
\frac{p}{1-(1/2)^{(2p)/(2-p)}},
\]
and the latter expression is bounded uniformly in $0<p\leq q$.
\end{pf}

In particular, the result reveals that the LHS of
(\ref{eq:nonconv}) is bounded by
%
\begin{equation}\label{eq: bound1}
\hat{d}_{2,N}(\hat{A},A^*)^{1-p/(2-p)}
\Vert\hat{A}-A^*\Vert_{S_p}^{p/(2-p)}
\sqrt{\vartheta/p}\biggl(\frac{M}{N}\biggr)^{1/2}
\end{equation}
with probability at least $1- C\exp(-\vartheta M/C^2)$ for any
$\sqrt{\vartheta}\geq C(q)$.

We now use the following simple consequence of the concavity of the
logarithm which is stated, for instance, in Tsybakov and van de Geer
(\citeyear{tsyfvdgeerff05}) (Lemma 5).

\begin{lemma}\label{lemma: d}
For any positive $v$, $t$ and any $\kappa
\geq1$, $\delta>0$ we have
\[
vt^{1/(2\kappa)}\leq(\delta/2)t +
c_{\kappa}\delta^{-1/(2\kappa-1)}v^{2\kappa/(2\kappa-1)},
\]
where $c_{\kappa}=(2\kappa-1)(2\kappa)\kappa^{-1/(2\kappa-1)}$.
\end{lemma}

Taking in Lemma \ref{lemma: d}
\[
t= \hat{d}_{2,N}(\hat{A},A^*)^2,\qquad v=
\Vert\hat{A}-A^*\Vert_{S_p}^{p/(2-p)}\sqrt
{\vartheta/p}\biggl(\frac{M}{N}\biggr)^{1/2},
\]
and $\kappa=(2-p)/(2-2p)$ shows that for any $\delta>0$
\[
(\ref{eq: bound1})\leq(\delta/2)\hat{d}_{2,N}(\hat{A}, A^*)^2+
\tau_7\delta^{p-1}\Vert\hat{A}-A^*\Vert_{S_p}^{p}
\]
with probability at least $1- C\exp(-\vartheta M/C^2)$.

The case $m\not= T$ can be deduced from the above result by the
following observation. For any matrix $B=(b_{ij})\in\R^{m\times T}$,
define the extension $\tilde{B}=(\tilde{b}_{ij})\in\R^{M\times M}$ with
$M=\max(m,T)$ as follows: $ \tilde{b}_{ij} = b_{ij}$ for $1\leq i\leq
m$, $1\leq j\leq T$ and $\tilde{b}_{ij} =0$ otherwise. Then one easily
checks that $ \Vert\tilde{B}\Vert_{S_p}=\Vert
B\Vert_{S_p}$ for all $p\in[0,\infty]$. Furthermore,
$\trace(B'X_i)=\trace(\tilde{B}'\tilde{X}_i)$ and
\[
\sup_{A\in\R^{M\times M}\setminus
\{0\}}\frac{N^{-1}\sum_{k=1}^N\trace(\tilde{X}_i'A)^2}{\Vert
A\Vert_{S_2}^2}=\sup_{A\in\R^{m\times
T}\setminus\{0\}}\frac{\Vert A\Vert_{2,N}^2}{\Vert
A\Vert_{S_2}^2} \leq c_0.
\]
Consequently, the result follows now from the already established proof
for the case $m=T$.
\end{pf*}


\section{Entropy numbers for quasi-convex Schatten class embeddings}
\label{sec: embeddings}

Here we derive bounds for the $k$th entropy numbers of the embeddings
$S_p^M\hookrightarrow S_2^M$ for $0<p<1$, where $S_p^M$ denotes the
$p$th Schatten class of real $M\times M$-matrices. Corresponding
results for the $l_p^M \hookrightarrow l_2^M$-embeddings are given
first by Edmunds and Triebel (\citeyear{edmfftriff89}) but their proof does not carry over
to the Schatten spaces. Pajor (\citeyear{pajorf98}) provides bounds for the
$S_p^M\hookrightarrow S_2^M$-embeddings in the convex case, $p\geq1$.
His approach is based on the trace duality (H\"{o}lder's inequality for
$p^{-1}+q^{-1}=1$) and the geometric formulation of Sudakov's
minoration
\[
\varepsilon\sqrt{\log\NN(A,\vert
\cdot\vert_2,\varepsilon)}\leq c\E\sup_{t\in A}\langle
G,t\rangle
\]
for some positive constant $c$, with a $d$-dimensional standard
Gaussian vector $G$ and an arbitrary subset $A$ of $\R^d$. Here
$\vert\cdot\vert_2 $ is the Euclidean norm in $\R^d$ and
$\langle\cdot,\cdot\rangle$ is the corresponding scalar product.
Gu\'{e}don and Litvak (\citeyear{glff00}) derive a slightly sharper
bound for the $l_p\hookrightarrow l_q$-embeddings than Edmunds and
Triebel (\citeyear{edmfftriff89}) with a different technique. In addition, they prove
lower bounds. We adjust their ideas concerning finite $\ell_p$ spaces
to the nonconvex Schatten spaces.

We denote by $e_k(\mathit{id}_{p,r}^M)$ the $k$th entropy number of the
embedding $S_p^M\hookrightarrow S_r^M$ for $0<p<r\leq\infty$, that is,
the infimum of all $\varepsilon>0$ such that there exist $2^{k-1}$
balls in $S_r^M$ of radius $\varepsilon$ that cover $\BB(S_p^M)$. For
the general definition of $k$th entropy numbers $e_k(T\dvtx F\rightarrow E)$
of bounded linear operators $T$ between quasi-Banach spaces $F$ and
$E$, we refer to Edmunds and Triebel (\citeyear{edmfftriff96}).

Recall that a homogeneous nonnegative functional $\Vert
\cdot\Vert$ is called $C$-quasi-norm, if it satisfies for all
$x,y$ the inequality $\Vert x+y\Vert\leq C\max(\Vert
x\Vert,\Vert y\Vert)$. Finally, any $p$-norm is a
$C$-quasi-norm with $C=2^{1/p}$ [cf., e.g., Edmunds and Triebel (\citeyear{edmfftriff96}),
page 2]. We will use the following lemma.

\begin{lemma}[[Gu\'{e}don and Litvak (\citeyear{glff00}\mbox{)]}]\label{lemma7}
Assume that $\Vert\cdot\Vert_i$ are symmetric
$C_i$-quasi-norms on $\R^n$ for $i=0,1$, and for some $\theta\in(0,1)$,
$\Vert\cdot\Vert_\theta$ is a quasi-norm on $\R^n$ such that
$\Vert x\Vert_{\theta}\leq\Vert
x\Vert_0^{\theta}\Vert x\Vert_1^{1-\theta}$ for all
$x\in\R^n$. Then for any quasi-normed space $F$, any linear operator
$T\dvtx F\rightarrow\R^n$, and all integers $k$ and $m$, we have
\[
e_{m+k-1}(T\dvtx F\rightarrow E_{\theta})\leq\bigl(C_0e_m(T\dvtx F\rightarrow
E_0)\bigr)^{\theta}\bigl(C_1e_k(T\dvtx F\rightarrow
E_1)\bigr)^{1-\theta},
\]
where $E_t$ stands for $\R^n$ equipped with quasi-norm $\Vert
\cdot\Vert_t$, $t\in\{0,\theta,1\}$.
\end{lemma}

Gu\'{e}don and Litvak (\citeyear{glff00}) did not specify the notion of symmetry
they used. So we have to clarify that here a (quasi-)norm $\Vert
\cdot\Vert$ is called symmetric if $(\R^n,\Vert
\cdot\Vert)$ is isometrically isomorphic to a symmetrically
(quasi-)normed operator ideal. This includes the diagonal operator
spaces (finite $\ell_p$) as well as the Schatten spaces. The proof of
Lemma \ref{lemma7} follows the lines of Pietsch (\citeyear{pietsch80}), Proposition 12.1.12,
replacing the triangle inequality by the quasi-triangle inequality.
Recall that the Schatten classes $S_p$ form interpolation couples like
their commutative analogs $\ell_p$.

\begin{lemma}[(Interpolation inequality)]\label{lemma7a}
For $0<p<q<r< \infty$, let $\theta\in[0,1]$ be such that
\[
\frac{\theta}{p}+\frac{1-\theta}{r}= \frac{1}{q}.
\]
Then, for all $A\in\R^{m\times T}$,
\[
\Vert A\Vert_{S_q}\leq\Vert
A\Vert_{S_p}^{\theta}\Vert A\Vert_{S_r}^{1-\theta}.
\]
\end{lemma}

Proof is immediate in view of the inequalities
\[
\sum_j a_j^q = \sum_j a_j^{\theta q} a_j^{(1-\theta) q} \le\biggl(\sum_j
a_j^p\biggr)^{\theta q/p} \biggl(\sum_j a_j^r\biggr)^{(1-\theta) q/r}
\]
valid for any nonnegative $a_j$'s.

\begin{proposition}[(Entropy numbers)]\label{prop: entropy}
Let $0<p<1$, $p<r\leq\infty$. Then there exists an absolute constant
$\beta$ independent of $p$ and $r$, such that for all integers $k$ and
$M$ we have
\[
e_k(\mathit{id}_{p,r}^M)\leq\min\biggl\{1, \alpha(\beta,p,r)\biggl(\frac
{M}{k}\biggr)^{1/p-1/r}\biggr\}
\]
with
\[
\alpha(\beta,p,r)\leq
2^{1+1/r}\biggl(\frac{\beta}{p}\biggr)^{1/p-1/r}\biggl(\frac{1}{1-p}\biggr)^{(1/p-1)(1/p-1/r)}.
\]
\end{proposition}


\begin{pf}
The fact that $e_k(\mathit{id}_{p,r}^M)$ is
bounded by $1$ is obvious, since $\BB(S_p^M)\subset\BB(S_r^M)$.
Consider the other case. We start with $r=\infty$ and then extend the
result to $r<\infty$ by interpolation. Fix some number $L> M$ and let
$D=D(M,L,p)$ be the smallest constant which satisfies, for all $1\leq
k\leq L$,
%
\begin{equation}\label{eq: ent4}
e_k(\mathit{id}_{p,\infty}^M) \leq D\biggl(\frac{M}{k}\biggr)^{1/p}.
\end{equation}
Let us show that $\alpha=\sup_{M,L}D(M,L,p)$ is finite. Since
$\Vert\cdot\Vert_{S_p}$, $p<1$, can be viewed as a
quasi-norm on $\R^{M^2}$ (isomorphic to $\R^{M\times M}$), Lemma
\ref{lemma7} applies with $F=E_0=S_p^M$, $E_1=S_{\infty}^M$,
$\theta=p$, $E_{\theta}=S_1^M$ and $m=1$. This gives
%
\begin{equation}
e_k(\mathit{id}_{p,1}^M)\leq4 (e_k(\mathit{id}_{p,\infty}^M))^{1-p}.\label{eq: ek1}
\end{equation}
Here the factor 4 follows from the relations $C_1=2$ and $C_0^p \leq
2$. Now, (\ref{eq: ek1}) and the factorization theorem for entropy
numbers of bounded linear operators between quasi-Banach spaces [see,
e.g., Edmunds and Triebel (\citeyear{edmfftriff96}), page 8], with factorization via
$S_1^M$, leads to the bound
%
\begin{eqnarray}\label{eq: ent3}
e_k(\mathit{id}_{p,\infty}^M)&\leq&
e_{[(1-p)k]}(\mathit{id}_{p,1}^M)e_{[pk]}(\mathit{id}_{1,\infty}^M)\nonumber
\\[-8pt]\\[-8pt]
&\leq&4\bigl(e_{[(1-p)k]}(\mathit{id}_{p,\infty}^M)\bigr)^{1-p}e_{[pk]}(\mathit{id}_{1,\infty
}^M),\nonumber
\end{eqnarray}
where for any $x\in(0,\infty)$, $[x]$ denotes the smallest integer
which is larger or equal to $x$.
Proposition 5 of Pajor (\citeyear{pajorf98}) 
entails $ \log\NN(\BB(S_1^M),\|\cdot\|_{S_{\infty}},\varepsilon
)\leq
cM/\varepsilon,\forall \varepsilon>0, $ and hence
%
\begin{equation}\label{eq: en2}
e_k(\mathit{id}_{1,\infty}^M)\leq c'M/k
\end{equation}
with constants $c$ and $c'$ independent of $M$, $\varepsilon$ and $k$.
Note that, in contrast to the $l_1^M\hookrightarrow
l_{\infty}^M$-embedding, for which the $k$th entropy numbers are
bounded by $ c''k^{-1}\log(1+M/k) $ with some $c''>0$ and $\log
_2M\leq
k\leq M$ [see, e.g., Edmunds and Triebel (\citeyear{edmfftriff96}), page 98], we have in
(\ref{eq: en2}) not a logarithmic but linear dependence of $M$ in the
upper bound. Plugging (\ref{eq: ent4}) and (\ref{eq: en2}) into
(\ref{eq: ent3}) yields
\begin{eqnarray*}
e_k(\mathit{id}_{p,\infty}^M)&\leq&4
\biggl(D\biggl(\frac{M}{(1-p)k}\biggr)^{1/p}\biggr)^{1-p} \frac{c'M}{pk}\\
&=&\frac{4c'}{p}\biggl(\frac{1}{1-p}\biggr)^{(1-p)/p}D^{1-p}\biggl(\frac{M}{k}\biggr)^{1/p}.
\end{eqnarray*}
Thus, by definition of $D$,
\[
D^p\leq\frac{4c'}{p}\biggl(\frac{1}{1-p}\biggr)^{(1-p)/p},
\]
which shows that $D$ is uniformly bounded in $M$ and $L$. This proves
the proposition for $r=\infty$.

Consider now the case $r<\infty$. In view of Lemma \ref{lemma7a} with
$\theta=p/r$, we can apply Lemma \ref{lemma7} with $F=E_0=S_p^M$,
$E_1=S_{\infty}^M$, $\theta=p/r$, $E_{\theta}=S_r^M$ and $m=1$. This
yields
\begin{eqnarray*}
e_k(\mathit{id}_{p,r}^M)&\leq&2^{1+1/r}(e_k(\mathit{id}_{p,\infty}^M))^{1-p/r}
\\
&\leq&2^{1+1/r}D^{1-p/r}\biggl(\frac{M}{k}\biggr)^{1/p-1/r}.
\end{eqnarray*}
\upqed\end{pf}

\begin{corollary}\label{cor: entropy}
For any $p\in(0,1)$, there exists a positive constant $\alpha_0(p)$
such that for all integers $M\geq1$ and any $\varepsilon\in(0,1]$,
\[
\log\NN(\BB(S_p^{M}),\|\cdot\|_{S_2},\varepsilon)\leq\alpha_0(p)
M\varepsilon^{-2p/(2-p)}.
\]
Moreover, $\alpha_0(p)=\mathrm{O}(1/p)$ for $p\searrow0$.
\end{corollary}

\begin{pf}
The result follows by transforming the
entropy number bound of Proposition \ref{prop: entropy} into an entropy
bound. Specification of the constant in Proposition \ref{prop: entropy}
yields
\[
\alpha_0(p)=\mathrm{O}\biggl( \frac{\beta}{p}\biggl(1+\frac{p}{1-p}\biggr)^{(1-p)/p}\biggr)=\mathrm{O}(1/p)
\]
as $p\searrow0$.
\end{pf}

\section*{Acknowledgment}
We are grateful to Alain Pajor for
pointing out reference Gu\'{e}don and Litvak (\citeyear{glff00}). 


%

\printaddresses

\end{document}